\documentclass[10pt]{article}
\usepackage[T1]{fontenc}
\usepackage{amsmath,amsfonts,amsthm,mathrsfs,amssymb}
\usepackage{cite}
\usepackage{graphicx,float}
\usepackage{subcaption}
\usepackage{graphicx}
\usepackage{placeins}
\usepackage{xcolor}
\usepackage{indentfirst}
\usepackage[bookmarks=true]{hyperref}
\usepackage{bm}
\numberwithin{equation}{section}
\newcommand{\R}{{\mathbb R}}

\newcommand{\N}{{\mathbb N}}
\newcommand{\C}{{\mathbb C}}

\newcommand{\bb}{\mathbb}
\newcommand{\bs}{\boldsymbol}

\newcommand{\be}{\begin{eqnarray}}
\newcommand{\ben}{\begin{eqnarray*}}
\newcommand{\en}{\end{eqnarray}}
\newcommand{\enn}{\end{eqnarray*}}

\newcommand{\pa}{\partial}

\newcommand{\ov}{\overline}
\newcommand{\curl}{{\rm curl\,}}

\newcommand{\divv}{{\rm div\,}}

\newcommand{\br}{\mathbf{r}}

\newcommand{\bpsi}{\bm{\psi}}
\newcommand{\bphi}{\bm{\phi}}

\newcommand{\bx}{\bm{x}}
\newcommand{\by}{\bm{y}}
\newcommand{\bnu}{\bm{\nu}}

\newcommand{\half}{\frac{1}{2}}

\newtheorem{theorem}{Theorem}[section]
\newtheorem{lemma}[theorem]{Lemma}
\newtheorem{corollary}[theorem]{Corollary}

\newtheorem{remark}[theorem]{Remark}

\newtheorem{problem}[theorem]{Problem}

\definecolor{rot}{rgb}{0.000,0.000,0.000}

\begin{document}
\renewcommand{\theequation}{\arabic{section}.\arabic{equation}}
\begin{titlepage}
\title{Spectral Galerkin method for solving elastic wave scattering problems with multiple open arcs}

\author{Carlos Jerez-Hanckes\thanks{Facultad de Ingenier\'{i}a y Ciencias, Universidad Adolfo Ib\'a\~{n}ez, Santiago, Chile. Email: {\tt carlos.jerez@uai.cl}}\;,
Jos\'{e} Pinto\thanks{Facultad de Ingenier\'{i}a y Ciencias, Universidad Adolfo Ib\'a\~{n}ez, Santiago, Chile. Email: {\tt jose.pinto@uai.cl}}\;,
Tao Yin\thanks{LSEC, Institute of Computational Mathematics and Scientific/Engineering Computing, Academy of Mathematics and Systems Science, Chinese Academy of Sciences, Beijing 100190, China. Email: {\tt yintao@lsec.cc.ac.cn}}}
\end{titlepage}
\maketitle

\begin{abstract}
We study the elastic time-harmonic wave scattering problems on unbounded domains with boundaries composed of finite collections of disjoints finite open arcs (or cracks) in two dimensions. Specifically, we present a fast spectral Galerkin method for solving the associated weakly- and hyper-singular boundary integral equations (BIEs) arising from Dirichlet and Neumann boundary conditions, respectively. Discretization bases of the resulting BIEs employ weighted Chebyshev polynomials that capture the solutions' edge behavior. We show that these bases guarantee exponential convergence in the polynomial degree when assuming analyticity of sources and arcs geometries. Numerical examples demonstrate the accuracy and robustness of the proposed method with respect to number of arcs and wavenumber.
\end{abstract}

\section{Introduction}
\label{sec:1}

We study the elastic time-harmonic wave scattering problems on unbounded domains with boundaries composed of finite collections of disjoints finite open arcs (or cracks) in two dimensions. Such problems play fundamental roles in multiple important applications in science and engineering such as in non-destructive testing of solid materials; detection of fractures; energy production from natural gas and geothermal resources \cite{MG14,PG15,W06}; to name a few. Thus, developing fast, accurate and efficient numerical schemes that can deal simultaneously with large numbers of arcs and a broad range of wavelengths is of particular interest for these applications. Among many available choices, we will focus on boundary integral equation (BIE) methods as they only require discretization on the domain boundaries and enforce the radiation condition at infinity automatically.

In this paper, we propose a spectral Galerkin method for solving both weakly- and hyper-singular BIEs resulting from two-dimensional elastic problems on multiple open arcs with Dirichlet and Neumann boundary conditions, correspondingly. We prove exponential convergence of the method by carrying out a comprehensive study of the spectral convergence analysis when both boundary data and arc geometries are given by analytic functions.

Numerical schemes for BIEs of open arc problems have been extensively studied for Laplace/Helmholtz
\cite{AS91,BL15,JP20,JP22,SW84,WS90}, 
elastostatic/elastodynamic \cite{CKM00,BXY21,HSW91,WS90} and Maxwell equations \cite{HU202}. Generally, their study requires handling the following three groups of questions:
\begin{itemize}
\item[(i)]  Are the BIEs well posed?
\item[(ii)] Is the numerical discretization consistent? Does it converge? How fast?
\item[(iii)] Is the associated linear system ill-conditioned? Is there need for preconditioning or regularization?
\end{itemize}

Under our setting of interest---multiple-arcs elastic wave scattering problems---, to answer (i) we will extend the approach employed for studying single-arc problems~\cite{HSW91,WS90}. {More precisely, we will} show that volume solutions can be expressed as superpositions of single- and double-layer potentials applied to surface densities over each arc~\cite{JP20} for Dirichlet and Neumann boundary conditions, respectively. Then, the corresponding weakly- and hyper-singular BIEs are derived after taking traces of these unknown potentials. Wellposedness of single-arc problems can then be extended to the multiple-arcs case by means of the Fredholm alternative (see Theorem~\ref{wellposedness} for the here considered case).

With respect to (iii), it should be pointed out that the resulting weakly- and hyper-singular BIEs are all of first kind, and thus, employing standard discretization bases leads to poor performance of iterative solvers for the linear systems arising from large scale problems (cf.~\cite{BL15} and references within). Then, one requires suitable preconditioners or regularized BIEs to construct more efficient numerical solvers. This topic has received attention in recent years~\cite{BL12,HJU14,HU202,EIJ21} and some attempts have been carried out to tackle elastic wave problems~\cite{BXY21,BY20,CDL21,XY22} by considering the composition of the weakly- and hyper-singular boundary integral operators (BIOs). Indeed, including
 preconditioning techniques into the presented numerical method is relevant but for the sake for brevity will be left as future work.

The present work mainly focus on tackling the second issue (ii) for the multiple arcs elastic problems though (i) is fully addressed for completeness . In this context, Nystr\"om-type strategies~\cite{CKM00} and  variational methods such as the boundary element method (BEM)~\cite{JP20} are traditionally employed for the numerical approximations of resulting open-arc BIEs. The Nystr\"om-type method has been developed for the two-dimensional elastostatic hyper-singular open-arc BIE in~\cite{CKM00} together with a convergence analysis in H\"older spaces. Although it is remarked in~\cite{CKM00} that the exponentially convergence holds true for the case of analytic arcs, the  square-root singularities at arc endpoints~\cite{CDD03}of the solutions are not considered. In fact, the edge singularities are considered recently in~\cite{BXY21} for the elastic open-arc BIEs, though no convergence estimates is provided. These singularities also hinder the performance of standard low-order BEM. Specifically, only suboptimal convergence can be obtained by using low-order uniform-mesh discretizations and additional techniques---graded or adaptive mesh refinement~\cite{FFHKP,VS90}, approximation space augmenation~\cite{SW84}---are required for improved convergence rates. Inspired by the spectral Galerkin-Bubnov discretization method proposed in \cite{AS91} for logarithmic kernel singular BIEs on a single-arc, a novel spectral Galerkin method is recently developed in \cite{JP20} for the numerical discretization of weakly-singular BIEs for Laplace and Helmholtz multiple finite arc problems. Therein, the approximation basis is given by weighted first-kind Chebyshev polynomials and rigorous error convergence estimates are proven based on the asymptotic properties of the Fourier-Chebyshev expansions. This leads to exponential convergence rates when both arcs and sources can be represented by analytic functions.

In this work, we extend the spectral Galerkin method \cite{JP20} to the more challenging case of elastic wave scattering by multiple open-arcs. Analogous to~\cite{CKM00}, an adequate Maue's representation formula~\cite{BXY21,YHX17} for elastodynamic problems is used to simplify the discretization of the hyper-singular BIE. Yet, and unlike~\cite{CKM00}, the corresponding variational formulation of the hyper-singular BIE avoids the treatment of tangential derivatives of weakly-singular operators. Due to the diverse edge singularities of the solutions of the weakly- and hyper-singular BIEs and so as to avoid any arc meshing, weighted first- and second-kind Chebyshev polynomials are utilized to construct the approximation basis, respectively. Then, by examining the polynomial expansion of the BIEs solutions, rigorous exponential convergence in the polynomial degree is proven by assuming the analyticity of the open-arcs and sources. In contrast to the traditional convergence analysis for BEM on arcs problems, we do not rely on localizing solution singularities by means of smooth window functions. Hence, we are able to obtain  exponential convergence as the window function is not analytic. This convergence analysis substantively improves the analysis of open-arc BIEs and we believe that it can provide a new strategy to prove convergence for the corresponding Nystr\"om-type methods while, unlike~\cite{CKM00}, the edge singularities are explicitly involved in the approximation spaces~\cite{BL12,BXY21}.

The remainder of this paper is organized as follows. Sections~\ref{sec:geo} through \ref{sec:traces} set forward formal definitions while Section~\ref{sec:elastic} describes the elastic open-arcs scattering problems here considered along with their corresponding weakly- and hyper-singular BIEs and their wellposedness. Approximation spaces and reduced linear systems of the spectral Galerkin numerical scheme are introduced in Section~\ref{sec:3.1}. We prove the exponential convergence of the method in Section~\ref{sec:3.3} through the investigation of the Chebyshev regularity of the solutions of BIEs (see Section~\ref{sec:3.2}). Numerical experiments illustrating the accuracy of the method are presented in Section~\ref{sec:4} including implementation details---compression algorithm---for multiple arcs problems.

\section{Definitions and problem statement}
\label{sec:definitions}

Vectors will be denoted in bold face, e.g., $\bs u$, and their components as $\bs u =(u^1,u^2)^\top$. We will say a function in $[-1,1]$ is analytic if it has an analytic complex extension to an open neighborhood of $[-1,1] \subset \mathbb{C}$. In particular, this implies that the function needs to have a complex extension to a Bernstein ellipse\footnote{The ellipse in the complex plane with foci $\pm1$, and major and minor semi-axes $\frac{1}{2}(\varrho + \varrho^{-1})$ and  $\frac{1}{2}(\varrho - \varrho^{-1})$, respectively.} of parameter $\varrho$, for some $\varrho >1$.

\subsection{Geometry}
\label{sec:geo}

We define open arcs (cracks) as injective functions $\bs r : [-1,1] \rightarrow \mathbb{R}^2$, such that each component is continuously differentiable, and also $\| \bs r'(t)\| >0$ for every $t \in [-1,1]$. Slightly abusing notation, we also call open arc the range of a function with the properties described above and the corresponding function is referred as arc parametrization. Thus, for any open arc a parametrization is implicitly fixed. Notice that under this convention two arcs are equal if their parametrization are the same, and not if only if their corresponding ranges coincide. Furthermore, we will say that an open arc is analytic if both of its associated parametrization coordinates are analytic.

Throughout $M$ is be a fixed natural number and $\Gamma := \{ \Gamma_1, \hdots, \Gamma_M\}$ denotes a set of $M$ disjoint analytic open arcs $\Gamma_i$, with parametrizations denoted by $\bm r_i$, for $i =1, \hdots,M.$ We assume that for each open arc $\Gamma_i$ there exists closed arc such that $\Gamma_i\subset\widetilde{\Gamma}_i$ with $\widetilde{\Gamma}_i$ being the boundary of a bounded subdomain $\Omega_i$ for which an orientation exists.
\subsection{Sobolev spaces}
\label{sec:sobolev}
We recall the standard Sobolev framework for open arcs~\cite[Section 2.3]{JP20}. For $G\subseteq\R^d$, $d=1,2$, being an open domain, and $s\in\R$, we denote by $H^s(G)$ the standard Sobolev spaces in $L^2(G)$ and by $H_{\rm loc}^s(G)$ their locally integrable counterparts. For an open arc $\Lambda$, we assume that there exists a closed arc $\widetilde{\Lambda}$ that contains $\Lambda$ and denote by $H^s(\widetilde{\Lambda})$ the Sobolev spaces defined through local parametrizations. We further define
\begin{align*}
H^s(\Lambda)&:= \{ u \in \mathcal{D}^*(\Lambda): \exists \ U  \in H^s(\widetilde{\Lambda}), u = U|_\Lambda \}, \\
\widetilde{H}^s(\Lambda)&:=\{u\in H^s(\widetilde{\Lambda}): \textrm{supp}(u) \subset \overline{\Lambda}\}.
\end{align*}
Here, $\mathcal{D}^*(\Lambda)$ denotes the dual space---with respect to the dual product $\langle \cdot, \cdot \rangle_{\Lambda}$---to $C^\infty_0(\Lambda)$, the space of infinitely differentiable functions with compact support on $\Lambda$. One can identify Sobolev dual spaces as follows
\ben
\widetilde{H}^{-s}(\Lambda)=(H^s(\Lambda))^*,\quad H^{-s}(\Lambda)=(\widetilde{H}^s(\Lambda))^*.
\enn
For the finite union of disjoint open arcs $\Gamma$, we define piecewise spaces as
\ben
\mathbb{H}^{s}(\Gamma):=\prod_{i=1}^M H^s(\Gamma_i)\times H^s(\Gamma_i),
\enn
and similarly for spaces $\widetilde{\mathbb{H}}^s(\Gamma)$, for $s\in\R$. The duality between $\mathbb{H}^{s}(\Gamma)$, and $\widetilde{\mathbb{H}}^{-s}(\Gamma)$ is given by
\ben
\langle \bs u, \bs v \rangle_\Gamma = \sum_{i=1}^M \langle u_i^1 , v_i^1 \rangle_{\Gamma_i} +\langle u_i^2 , v_i^2 \rangle_{\Gamma_i}.
\enn
\subsection{Elasticity Dirichlet and Neumann traces}
\label{sec:traces}
Let us introduce the Dirichlet and Neumann traces for elastodynamics following \cite{JP20,M00}. For an open arc $\Gamma_i \in \Gamma$, and a $C^\infty$-function $U$ that is smooth on a neighborhood of $\Gamma_i$, we define the interior $(-)$ (resp.~exterior $(+)$) Dirichlet traces:
\ben
\gamma_{D,i}^\pm U(\bx):=\lim_{\varepsilon\rightarrow 0} U(\bx \pm \varepsilon \bnu_i)\quad \forall\ \bx\in\Gamma_i,
\enn
where $\bnu_i$ denotes the unitary normal vector with direction $(r_{i,2}',-r_{i,1}')^\top$. If $\gamma_{D,i}^+U=\gamma_{D,i}^-U$, we write $\gamma_{D,i}U=\gamma_{D,i}^\pm U$. These definitions can be extended to more general Sobolev spaces by density. In particular, we know that $\gamma_{D,i}^\pm:H_{\rm loc}^1(\Omega)\rightarrow H^{\half}(\Gamma_i)$ is bounded \cite[Theorem 3.37]{M00}.

For a smooth vector ${\bm U}=(U^1,U^2)^\top$, its Dirichlet trace is given by those of its components, i.e.~$\gamma_{D,i}^\pm {\bm U}=(\gamma_{D,i}^\pm U^1,\gamma_{D,i}^\pm U^2)^\top$ and thus, $\gamma_{D,i}^\pm:H_{\rm loc}^1(\Omega)^2\rightarrow H^{\half}(\Gamma_i)^2$ is also bounded. In contrast to the Laplace or Helmholtz cases~\cite{JP20},  elasticity Neumann traces are defined in terms of the traction operator $\mathcal{T}(\pa,\boldsymbol{\nu})$ defined by
\be
\label{traction}
\mathcal{T}(\pa,\boldsymbol{\nu}){\bm U}:=2 \mu \, \partial_{\bnu} {\bm U} + \lambda \,
\bnu \, \divv {\bm U}-\mu \bnu^\perp\curl {\bm U},
\en
in which $\lambda, \mu$ are Lam\'e parameters, $\bnu^\perp:=(-\nu^2,\nu^1)^\top$, $\partial_{\bnu}:=\bnu\cdot\nabla$ is the normal derivative and the two-dimensional scalar operator $\curl$ is defined as $\curl\bm{U}:=\pa_1U^2-\pa_2U^1$. Then, we can define Neumann traces for smooth vector fields ${\bm U}$ as
\ben
\gamma_{N,i}^\pm\bm{U}:= \lim_{\varepsilon\rightarrow 0} \mathcal{T}(\pa,\bm \nu_i){\bm U}(\bx\pm \varepsilon\bnu_i)\quad \forall \ \bx\in\Gamma_i.
\enn
The Neumann trace can be extended to a bounded map: $\gamma_{N,i}^\pm:H_{\rm loc}^1(\Omega)^2 \cap  \{ {\bs U} : \Delta^* \bs U \in  L^2_{loc}(\Omega)^2 \} \rightarrow H^{-\half}(\Gamma_i)^2$, wherein  $\Delta^*$ is defined in \eqref{LameOper}, \cite[Chapter 4]{M00}. As for the Dirichlet case, if $\gamma_{N,i}^+\bm U=\gamma_{N,i}^-\bm U$, we denote $\gamma_{N,i}\bm U=\gamma_{N,i}^\pm\bm U$.

\subsection{Elastic wave scattering problems and BIEs}
\label{sec:elastic}
We consider the problem of elastic time-harmonic wave scattering in the unbounded domain $\Omega:= \mathbb{R}^2\setminus \Gamma$ which is filled with a linear isotropic and homogeneous solid medium characterized by the Lam\'e constants $\lambda, \mu$ with $\mu>0$, $\lambda+\mu>0$, and mass density $\rho>0$. In particular, we seek the displacement field ${\bs U}=(U^1,U^2)^\top\in H_{\rm loc}^1(\Omega)^2$ satisfying the time-harmonic Navier equation
\be
\label{navier}
\Delta^*{\bm U}+\rho\omega^2\bm{U}={\bm 0} \quad\mbox{in}\quad\Omega,
\en
and the Kupradze-Sommerfeld radiation condition~\cite{KGBB79} at infinity where $\omega>0$ denotes the angular frequency. Here, $\Delta^{*}$ is the Lam\'e operator given by
\be
\label{LameOper}
\Delta^* := \mu\,\mbox{div}\,\mbox{grad} + (\lambda + \mu)\,\mbox{grad}\, \mbox{div}.
\en
On $\Gamma$ the solution is assumed to satisfy either the Dirichlet boundary condition
\be
\label{Dirichlet}
\gamma_{D,i}^\pm{\bm U}={\bm f}_i \quad\mbox{on}\quad\Gamma_i,\ i=1,\ldots, M,
\en
or the Neumann counterpart
\be
\label{Neumann}
\gamma_{N,i}^\pm{\bm U}={\bm g}_i \quad\mbox{on}\quad\Gamma_i,\ i=1,\ldots, M,
\en
for ${\bm f}_i$ and ${\bm g}_i$ in $H^\frac{1}{2}(\Gamma_i)$ and $H^{-\frac{1}{2}}(\Gamma_i)$, respectively. It follows from~\cite{HSW91,SW84,WS90} that the solutions of (\ref{navier}) under Dirichlet and Neumann boundary conditions can be expressed in terms of either single- and double-layer potentials over cracks $\Gamma_i$:
 \be
\label{DirichletS}
{\bs U}(\bx)  = \sum_{i=1}^M (\mathcal{S}_i\bphi_i)(\bx),\quad (\mathcal{S}_i\bphi_i)(\bx):=  \int_{\Gamma_i}\bb E(\bx,\by)\bphi_i(\by)\,ds_{\by}, \quad \forall\,\bx\in\Omega,
\en
and
\be
\label{NeumannD}
{\bs U}(\bx) = \sum_{i=1}^M (\mathcal{D}_i\bpsi_i)(\bx),\quad (\mathcal{D}_i\bpsi_i)(\bx):= \int_{\Gamma_i}(\mathcal{T}(\pa_{\by},\bnu_{\by})
{\mathbb E}(\bx,\by))^\top \bpsi_i(\by)\,ds_{\by}, \quad \forall\,\bx\in\Omega,
\en
respectively. Here, ${\mathbb E}(\bx,\by)$ denotes the fundamental displacement tensor for the Navier equation in $\R^2$ given by
\ben
{\mathbb E}(\bx,\by)=\frac{1}{\mu}\gamma_{\kappa_s}(\bx,\by){\mathbb I}+\frac{1}{\rho\omega^2}
\nabla_{\bx}\nabla_{\bx}^\top \left[\gamma_{\kappa_s}(\bx,\by)-\gamma_{\kappa_p}(\bx,\by)\right].
\enn
with $\bb I$ being the $2\times2$ identity, $\gamma_{k}(\bx,\by)$ denoting the fundamental solution of the Helmholtz equation in $\R^2$ with wavenumber $\kappa$,
\be
\label{HelmholtzFS}
\gamma_{k}(\bx,\by) = \frac{\iota}{4}H_0^{(1)}(\kappa |\bx-\by|), \quad \bx\ne \by.
\en
wherein $\iota$ denotes the imaginary unit, and  $H_0^{(1)}(\cdot)$ signals the zeroth-order Hankel function of first kind \cite[9.1.3]{stg}. The wavenumbers
\ben
\kappa_s := \omega\sqrt{\frac{\rho}{\mu}},\quad \kappa_p := \omega\sqrt{\frac{\rho}{\lambda+2\mu}}
\enn
correspond to elastodynamic compressional and shear waves, respectively. Moreover, the unknown densities $\bphi_i$ in (\ref{DirichletS}) and $\bpsi_i$ in (\ref{NeumannD}) represent the jumps of elastic Neumann and Dirichlet traces at $\Gamma_i$, respectively, i.e.
\ben
\bphi_i=-\gamma_{N,i}^+\bm{U}+\gamma_{N,i}^-\bm{U},\quad \bpsi_i=\gamma_{D,i}^+\bm{U}-\gamma_{D,i}^-\bm{U}.
\enn
Let us define the following boundary integral operators (BIOs), corresponding to weakly- and hyper-singular ones when $i=j$, as
\begin{align}
\label{SBIO}
\mathcal{V}_{ij}[\bphi_j](\bx)  &:={\gamma_{D,i}}_{\bx}\int_{\Gamma_j}{\bb E}(\bx,\by)\bphi_j(\by)ds_{\by},\quad \bx\in\Gamma_i,\\
\label{HBIO}
\mathcal{W}_{ij}[\bpsi_j](\bx) &:= \gamma_{{N,i}_{\bx}}\int_{\Gamma_j}
\gamma_{{N,i}_{\by}}
{\bb E}(\bx,\by))^\top\bpsi_j(\by)ds_{\by},\quad \bx\in\Gamma_i,
\end{align}
where the integral for the second operator is understood as a principal value. These BIOs are well defined regardless of the sign of the trace operation \cite[Chapter 6]{M00}. The original Dirichlet and Neumann volume problems can be reduced to the following BIEs on $\Gamma$:
\be
\label{DBIE}
\boldsymbol{\mathcal{V}}[\bphi]={\bs f}, \\
\label{NBIE}
\boldsymbol{\mathcal{W}}[\bpsi]={\bs g},
\en
respectively, wherein we have defined
\ben
\boldsymbol{\mathcal{V}}:=\begin{bmatrix}
\mathcal{V}_{11} & \mathcal{V}_{12} & \ldots & \mathcal{V}_{1M}\\
\mathcal{V}_{21} & \mathcal{V}_{22} & \ldots & \mathcal{V}_{2M}\\
\vdots & \vdots & \ddots & \vdots \\
\mathcal{V}_{M1} & \mathcal{V}_{M2} & \ldots & \mathcal{V}_{MM}
\end{bmatrix}\quad \boldsymbol{\mathcal{W}}:=\begin{bmatrix}
\mathcal{W}_{11} & \mathcal{W}_{12} & \ldots & \mathcal{W}_{1M}\\
\mathcal{W}_{21} & \mathcal{W}_{22} & \ldots & \mathcal{W}_{2M}\\
\vdots & \vdots & \ddots & \vdots \\
\mathcal{W}_{M1} & \mathcal{W}_{M2} & \ldots & \mathcal{W}_{MM}
\end{bmatrix},
\enn
and
\begin{align*}
\bphi&=(\bphi_1,\bphi_2,\ldots,\bphi_M)^\top,\quad \bpsi=(\bpsi_1,\bpsi_2,\ldots,\bpsi_M)^\top,\\
{\bs f}&=({\bs f}_1,{\bs f}_2,\ldots,{\bs f}_M)^\top,\quad {\bs g}=({\bs g}_1,{\bs g}_2,\ldots,{\bs g}_M)^\top.
\end{align*}
The boundary integral problems corresponding to the Dirichlet/Neumann elastic problems are summarized as follows:
\begin{problem}
\label{problem1}
Given ${\bs f}\in \mathbb{H}^{\half}(\Gamma)$ and ${\bs g}\in \mathbb{H}^{-\half}(\Gamma)$, we seek $\bphi\in\widetilde{\mathbb{H}}^{-\half}(\Gamma)$ and $\bpsi\in\widetilde{\mathbb{H}}^{\half}(\Gamma)$ such that
\ben
\boldsymbol{\mathcal{V}}[\bphi]={\bs f}, \quad \boldsymbol{\mathcal{W}}[\bpsi]={\bs g},
\enn
or equivalently,
\begin{align*}
\langle\boldsymbol{\mathcal{V}}[\bphi],{\bs v}\rangle_\Gamma&=\langle{\bs f},{\bs v}\rangle_\Gamma, \quad\forall \ {\bs v}\in \widetilde{\mathbb{H}}^{-\half}(\Gamma),\\
\langle\boldsymbol{\mathcal{W}}[\bpsi],{\bs v}\rangle_\Gamma&=\langle{\bs g},{\bs v}\rangle_\Gamma, \quad\forall \ {\bs v}\in \widetilde{\mathbb{H}}^{\half}(\Gamma).
\end{align*}
\end{problem}
The following lemma gives the coercivity of the diagonal operators $\mathcal{V}_{ii}$ and $\mathcal{W}_{ii}$ by transform the problem into a closed domain one~\cite{SW84}.
\begin{lemma}
\label{coercive}
There exist compact operators $\mathcal{C}_{ii}^1: \widetilde{H}^{-\half}(\Gamma_i)^2\rightarrow H^{\half}(\Gamma_i)^2$ and $\mathcal{C}_{ii}^2: \widetilde{H}^{\half}(\Gamma_i)^2\rightarrow H^{-\half}(\Gamma_i)^2$ such that for any $\bphi_i\in \widetilde{H}^{-\half}(\Gamma_i)^2$ and $\bpsi_i\in \widetilde{H}^{\half}(\Gamma_i)^2$, there exist positive constants $c_i^1,c_i^2$ such that
\begin{align*}
\left|\left\langle (\mathcal{V}_{ii}+\mathcal{C}_{ii}^1)[\bphi_i],\bphi_i\right\rangle_{\Gamma_i}\right| &\ge c_i^1\|\bphi_i\|_{\widetilde{H}^{-\half}(\Gamma_i)^2}^2,\\
\left|\left\langle (\mathcal{W}_{ii}+\mathcal{C}_{ii}^2)[\bpsi_i],\bpsi_i\right\rangle_{\Gamma_i}\right| &\ge c_i^2\|\bpsi_i\|_{\widetilde{H}^{\half}(\Gamma_i)^2}^2.
\end{align*}
\end{lemma}
\begin{proof}
Following \cite[Section 3.5.3]{SS10}, the coercivity properties for integral operators on open arcs can be deduced directly from the closed boundary case, and the latter is a well known result for the elastic wave operators, see  \cite[(6.5) and (6.10)]{HKR00} .
\end{proof}

Now we show that Problem \ref{problem1} is well posed.

\begin{theorem}
\label{wellposedness}
For any ${\bs f}\in \mathbb{H}^{\half}(\Gamma)$ and ${\bs g}\in \mathbb{H}^{-\half}(\Gamma)$, there exist unique solutions $\bphi\in\widetilde{\mathbb{H}}^{-\half}(\Gamma)$ and $\bpsi\in\widetilde{\mathbb{H}}^{\half}(\Gamma)$ for Problem \ref{problem1}. Moreover, the solution operators are bounded, i.e.
\ben
\|\bphi\|_{\widetilde{\mathbb{H}}^{-\half}(\Gamma)}\lesssim \|{\bs f}\|_{\mathbb{H}^{\half}(\Gamma)},\quad \|\bpsi\|_{\widetilde{\mathbb{H}}^{\half}(\Gamma)}\lesssim \|{\bs g}\|_{\mathbb{H}^{-\half}(\Gamma)}.
\enn
\end{theorem}
\begin{proof}
Following the same arguments as those for Lemma~\ref{coercive}, it can be easily verified that for $i,j\in\{1,\ldots,M\}$,
\begin{align*}
\mathcal{V}_{ij}: &\widetilde{H}^{-\half}(\Gamma_i)^2\rightarrow H^{\half}(\Gamma_j)^2,\\
\mathcal{W}_{ij}: &\widetilde{H}^{\half}(\Gamma_i)^2\rightarrow H^{-\half}(\Gamma_j)^2,
\end{align*}
are all bounded operators. In particular, if $i\ne j$, the operators $\mathcal{V}_{ij}$ and $\mathcal{W}_{ij}$ are compact as the kernel function is at least ${C}^1$ in each component. Thus, by the coercivity result of Lemma~\ref{coercive} and the Fredholm alternative~\cite[Theorem 2.33]{M00}, one only needs proving injectivity to ensure existence. For $M=1$, injectivity follows via the same arguments presented in \cite[Sections 2-3]{WS90}.

The injectivity for the general case of $M>1$ can be shown following the proof idea of \cite[Theorem 3.10]{JP20}. More precisely, let $\bphi=(\bphi_1,\bphi_2,\ldots,\bphi_M)^\top$ and $\bpsi=(\bpsi_1,\bpsi_2,\ldots,\bpsi_M)^\top$ be such that
\ben
\sum_{j=1}^M\mathcal{V}_{ij}[\bphi_j]=\bm{0},\quad \sum_{j=1}^M\mathcal{W}_{ij}[\bpsi_j]=\bm{0}\quad \forall \ i\in\{1,\ldots,M\},
\enn
and define the volume potentials
\ben
\bm{U}_j:=\mathcal{S}_j\bphi_j,\quad \bm{V}_j:=\mathcal{D}_j\bpsi_j.
\enn
These last ones are solutions of the elastic problems (\ref{navier}) over $\R^2\backslash\ov\Gamma_j$ as well as the superpositions $\bm{U}=\sum_{j=1}^M\bm{U}_j$ and $\bm{V}=\sum_{j=1}^M\bm{V}_j$ defined over $\Omega$. Then, it holds that
\ben
\gamma_{D,i}\bm{U}=0,\quad \gamma_{N,i}\bm{V}=0\quad \forall \ i\in\{1,\ldots,M\}.
\enn
The uniqueness of the elastic open-arc problems with zero Dirichlet/Neumann boundary condition implies that $\bm{U}=\bm{V}=\bm 0$ in $\Omega$, i.e.
\be
\label{proof2.3-1}
\bm{U}_i=-\sum_{j\ne i} \mathcal{S}_j\bphi_j,\quad \bm{V}_i=-\sum_{j\ne i} \mathcal{D}_j\bpsi_j.
\en
Now using the jump relations of the single and double layer potentials~\cite{HW08,KGBB79} we obtain that
\begin{align*}
\bphi_i&=  \gamma_{N,i}^-\mathcal{S}_{\Gamma_i}\bphi_i -\gamma_{N,i}^+\mathcal{S}_{\Gamma_i}\bphi_i =-\sum_{j\ne i} \left(\gamma_{N,i}^-\mathcal{S}_{\Gamma_j}\bphi_j -\gamma_{N,i}^+\mathcal{S}_{\Gamma_j}\bphi_j\right)=\bm 0,\\
\bpsi_i&= \gamma_{D,i}^+\mathcal{D}_{\Gamma_i}\bpsi_i -\gamma_{D,i}^-\mathcal{D}_{\Gamma_i}\bpsi_i  =-\sum_{j\ne i} \left(\gamma_{D,i}^+\mathcal{D}_{\Gamma_j}\bpsi_j -\gamma_{D,i}^-\mathcal{D}_{\Gamma_j}\bpsi_j\right)=\bm 0,
\end{align*}
where the right-most equalities follows since $\sum_{j\ne i} \mathcal{S}_j\bphi_j,$ and $\sum_{j\ne i} \mathcal{D}_j\bpsi_j$ are smooth functions, thus their jumps are zero.
\end{proof}

\section{Spectral Galerkin numerical scheme}
\label{sec:3}

We now describe a spectral Galerkin numerical scheme for solving Problem~\ref{problem1} and establish specific convergence rates extending our previous work for Laplace and Helmholtz problems \cite{JP20}.

\subsection{Approximation spaces}
\label{sec:3.1}

This section is devoted to constructing dense conforming high-order discretizations for the spaces $\widetilde{\mathbb{H}}^{\pm\half}(\Gamma)$. Specifically, we employ weighted Chebyshev polynomials per arc to generate high-order global polynomial bases. As is well known~\cite{CDD03}, the solutions of the weakly-singular and hyper-singular BIEs of Problem~\ref{problem1} admit square-root singularities at the arcs' endpoints and more precisely, the solutions $\bs\phi_i$ and $\bs\psi_i$ take the asymptotic behaviors as $\mathrm{dist}_i^{-1/2}$ and $\mathrm{dist}_i^{1/2}$ with, on each $\Gamma_i$, $\mathrm{dist}_i$ representing the distance to the endpoint of $\Gamma_i$. In particular, for the considered parameterization of the open-arcs, we take the scalar function $w(t):=\sqrt{1-t^2}(\sim \mathrm{dist}_i^{1/2})$, $t\in(-1,1)$ to reproduce the asymptotics of the solutions $\bs\phi_i$ and $\bs\psi_i$.

We denote by $\{T_n\}_{n=0}^N$ the set of first $N+1$ first-kind Chebyshev polynomials, orthogonal under the weight $w^{-1}$, and denote by $\{U_n\}_{n=0}^N$ the set of first $N+1$ second-kind Chebyshev polynomials, orthogonal under the weight $w$. It follows that
\begin{align}
\label{eq:chebortg}
\int_{-1}^1 T_n(t)T_m(t)w^{-1}(t)dt=\begin{cases}
0, & n\ne m,\cr
\pi, & n=m=0,\cr
\frac{\pi}{2}, & n=m\ne0,
\end{cases}\; \int_{-1}^1 U_n(t)U_m(t)w(t)dt=\begin{cases}
0, & n\ne m,\cr
\frac{\pi}{2}, & n=m.
\end{cases}
\end{align}

We consider the elements $p_n^i:={\dfrac{T_n\circ {\bs r}_i^{-1}}{|\br_i' \circ \br_i^{-1}|}}$ and $q_n^i:=U_n\circ {\bs r}_i^{-1}$ and the spaces they span are denoted by $\mathbf{ T}_N(\Gamma_i)$ and $\mathbf{ U}_N(\Gamma_i)$, respectively. We account for edge behavior by multiplying these bases by suitable weights and obtain the spaces
\begin{align*}
{\mathbf{ T}_N^w(\Gamma_i)}&:=\{\widetilde{p}^i=w_i^{-1}p^i: p^i\in {\bf T}_N(\Gamma_i)\},\\
{\mathbf{ U}_N^w(\Gamma_i)}&:=\{\widetilde{q}^i=w_iq^i: q^i\in {\bf U}_N(\Gamma_i)\},
\end{align*}
wherein $w_i=w\circ {\bs r}_i^{-1}$ and the corresponding bases for ${\bf T}_N^w(\Gamma_i)$ (resp.~${\bf U}_N^w(\Gamma_i)$) can be characterized as $\widetilde{p}_n^i=w_i^{-1}p_n^i$ (resp.~$\widetilde{q}_n^i=w_i q_n^i$).

For the case of multiple arcs, we define the following finite-dimensional approximation Cartesian product spaces:
\ben
\bb T_N:=\prod_{i=1}^M  \mathbf{T}_N^w(\Gamma_i)\times \mathbf{T}_N^w(\Gamma_i) ,\quad {\bb U}_N:=\prod_{i=1}^M \mathbf{U}_N^w(\Gamma_i)\times \mathbf{U}_N^w(\Gamma_i).
\enn

\begin{problem}[Linear system]
\label{problem2}
Let $M,N\in\N$. Given ${\bs f}\in \mathbb{H}^{\half}(\Gamma)$ and ${\bs g}\in \mathbb{H}^{-\half}(\Gamma)$, we seek coefficients ${\bs a}=(\bs a_1,\bs a_2,\ldots, \bs a_M)\in \C^{M2(N+1)}$ and  ${\bs b}=(\bs b_1, \bs b_2,\ldots,\bs b_M)\in \C^{M2(N+1)}$ such that
\be
\label{linearsys}
{\bs V}{\bs a}={\mathfrak{\bs f}},\quad {\bs W}{\bs b}={\mathfrak{\bs g}},
\en
wherein the entries of the matrix blocks ${\bs V}_{ij}\in \bb C^{{2}(N+1)\times {2}(N+1)}$ and ${\bs W}_{ij}\in \bb C^{{2}(N+1)\times {2}(N+1)}$ of the Galerkin matrices ${\bs V}\in \bb C^{M{2}(N+1)\times M{2}(N+1)}$ and ${\bs W}\in \bb C^{M2(N+1)\times M2(N+1)}$, respectively, are given by
\ben
({\bs V}_{ij})_{{lm,pq}}=\left\langle \mathcal{V}_{ij}[\widetilde{p}_m^j \bs e_p],\widetilde{p}_l^i\bs e_q\right\rangle_{\Gamma_i}\\
\quad ({\bs W}_{ij})_{{lm,pq}}=\left\langle \mathcal{W}_{ij}[\widetilde{q}_m^j \bs e_p],\widetilde{q}_l^i\bs e_q\right\rangle_{\Gamma_i}
\enn
for all $i,j=1,\ldots,M$; $l,m=0,\ldots,N$, and $p,q=1,2$, with $\bs e_1 = (1,0)^\top$, and $\bs e_2 =(0,1)^\top$. The corresponding discrete right-hand sides $\mathfrak{\bs f}=(\mathfrak{\bs f}_{1},\mathfrak{\bs f}_s,\ldots,\mathfrak{\bs f}_M)\in \C^{M2(N+1)}$ and $\mathfrak{\bs g}=(\mathfrak{\bs g}_{1},\mathfrak{\bs g}_s,\ldots,\mathfrak{\bs g}_M)\in \C^{M2(N+1)}$ have components $({\mathfrak{\bs f}}_{i})_{l,p}=\langle{\bs f}_i,\widetilde{p}_l^i\bs e_p\rangle_{\Gamma_i}$ and $(\mathfrak{{\bs g}}_{i}^N)_{l,p}=\langle{\bs g}_i,\widetilde{q}_l^i\bs e_p\rangle_{\Gamma_i}$ for all $i=1,\ldots,M$, $l=0,\ldots,N$, and $p=1,2$.
\end{problem}

By solving the linear systems (\ref{linearsys}) of Problem~\ref{problem2}, one can approximate solutions at each arc $\Gamma_i$ of  Problem~\ref{problem1} through the linear combinations:
\be
\label{eq:discretesols}
\bphi_i^N= \sum_{p=1}^2\sum_{l=0}^N (a_i)_{l,p}\widetilde{p}_l^i \bs e_p,\quad \bpsi_i^N= \sum_{p=1}^2 \sum_{l=0}^N (b_i)_{l,p}\widetilde{q}_l^i \bs e_p,
\en
for all $i=1,\ldots,M$. Denoting 
\ben
\bphi^N=(\bphi_1^N,\bphi_2^N,\ldots,\bphi_M^N)\quad\mathrm{and}\quad \bpsi^N=(\bpsi_1^N,\bpsi_2^N,\ldots,\bpsi_M^N), 
\enn
then the following quasi-optimality result holds for Galerkin discretizations (cf.~\cite[Section 4.2]{SS10}).

\begin{lemma}
\label{lemma:qopty}
There exists $N_0 \in \mathbb{N}$, such that for any $N>N_0$, the solutions $\bs a$, $\bs b$ of Problem \ref{problem2} exist, are unique, and the corresponding approximations $\bs \phi^N$ and $\bs \psi^N$ satisfy
\begin{align}
\label{quasi1}
\|{\bs \phi}-{\bs \phi}^N\|_{\widetilde{\mathbb{H}}^{-\half}(\Gamma)} &\lesssim \inf_{{\bs P}\in {\bb T}_N} \|{\bs \phi}-{\bs P}\|_{\widetilde{\mathbb{H}}^{-\half}(\Gamma)},\\
\label{quasi2}
\|{\bs \psi}-{\bs \psi}^N\|_{\widetilde{\mathbb{H}}^{\half}(\Gamma)} &\lesssim \inf_{{\bs Q}\in {\bb U}_N} \|{\bs \psi}-{\bs Q}\|_{\widetilde{\mathbb{H}}^{\half}(\Gamma)}.
\end{align}
\end{lemma}
While the quasi-optimality results \eqref{quasi1})-\eqref{quasi2} ensure convergence---provided that $N>N_0$---of the approximations, they do not provide any information on the speed of convergence. Indeed, in order to establish the rate of convergence of the spectral solver we expand the solutions $\bm{\phi}$, $\bm{\psi}$, of Problem \ref{problem1} as infinite series of the adequate Chebyshev polynomials. Then, we will show that the corresponding coefficients decay exponentially fast. In particular, we will consider the expansion on arcs
\ben
&& {\bs \phi}_i \circ \mathbf{r}_i(t) = \sum_{p=1}^2 \sum_{l=0}^\infty (a_i)_{l,p} w^{-1}(t)T_l(t){\bs e}_p , \quad i=1,\hdots,M,\\
&& {\bs \psi}_i \circ \mathbf{r}_i(t) = \sum_{p=1}^2 \sum_{l=0}^\infty (b_i)_{l,p} w(t)U_l(t){\bs e}_p, \quad i=1,\hdots,M,
\enn
and show that, under the assumption of analytic boundaries and right-hand-sides, it holds that
\begin{align*}
| (a_i)_{l,p}| \leq C \varrho^{-l}, \quad
| (b_i)_{l,p}| \leq C \varrho^{-l}, \quad \forall \ l \in \mathbb{N},\ p=1,2,\ i=1,\hdots,M,
\end{align*}
where $C, \varrho>1$ are generic constants that could be different for both equations, as well as for different values of $i$ and $p$. Moreover, these constants also depend on the geometry and problem parameters $\lambda, \mu, \rho, \omega$. With the decay rates of the coefficients  $(a_i)_{l,p}$ and $(b_i)_{l,p}$ at hand, the exponential convergence result
\ben
\|\bs \phi-\bs \phi^N\|_{\widetilde{\mathbb{H}}^{-\half}(\Gamma)} \leq C \varrho^{-N},\quad \|\bs \psi-\bs \psi^N\|_{\widetilde{\mathbb{H}}^{\half}(\Gamma)} \leq C \varrho^{-N}.
\enn
can be deduced from the quasi-optimality estimates \eqref{quasi1}--\eqref{quasi2}.

\begin{remark}
\label{laplace}
The idea on how to show the coefficients' asymptotic decay rates is based on a generalization of a trivial observation concerning the weakly-singular integral equation for the Laplace equation on the straight arc $\mathbf{r}(t) = (t,0)$, $t \in (-1,1)$. Specifically, consider the corresponding integral equation:
\be\label{eq:logint}
\int_{-1}^1 \log|s-t|\phi(t)dt =f(s),\quad s\in(-1,1),
\en
and note that by \cite[Lemma 4.5]{JP20}, one has
\be
\label{eq:logcofs}
\int_{-1}^1\int_{-1}^1 \log|s-t|w^{-1}(t)T_l(t) w^{-1}(s)T_n(s)dtds= {d_l :=}\begin{cases}
-\pi^2\log2, & l=n= 0, \cr
-\dfrac{\pi^2}{2n}, & l=n\ne 0, \cr
0, & l\ne n.
\end{cases}
\en
Then, if the solution $\phi$ is expanded as $\phi(t) = \sum_{n=0}^\infty a_l w^{-1}(t)T_l(t)$, the integral equation can be reduced to the following infinite system, with $a_l$ being the unknowns
\begin{align}
\label{eq:lapbie}
d_l a_l = f_l \quad l \in \mathbb{N},
\end{align}
where $|d_l|=\frac{\pi^2}{2l}, l\ge 1$ and $f_l  = \int_{-1}^1 f(t) w^{-1}(t)T_l(t)dt$ denotes the $l$th Chebyshev coefficient of the right-hand-side $f$. It is well known (cf.~\cite[Chapter 8]{T13}) that if $f$ is analytic, its coefficients decay exponentially, i.e.~$|f_l| \leq C \varrho^{-l}$, for some $\varrho > 1$. Hence, from \eqref{eq:lapbie} we deduce that
$$
|a_l| \leq C l \varrho^{-l},
$$
which can be expressed alternatively as $|a_l| \leq C \varrho^{-l},$ for a different $\varrho>1$.
\end{remark}

\subsection{Abstract Chebyshev regularity of solutions}
\label{sec:3.2}

The following result generalizes our previous remark for abstract weakly-singular integral equations.

\begin{lemma}
\label{lemma:coefsgeneralbies1}
Let $G(s,t)$ be a weakly-singular kernel which can be decomposed as
\begin{align}
\label{eq:kersplit1}
G(s,t) =  \log|s-t|J_0 + \log|s-t| (J(s,t)-J_0)+ R(s,t),\quad s\ne t,
\end{align}
with $J,R$ being analytic in both variables, $J_0 \in \mathbb{C}$. Additionally, assume that
$$
J(s,t) - J_0 = (s-t)^2 A(s,t),
$$
with $A$ also analytic in both variables. For an analytic function $f$, if the following integral equation
\begin{align}
\label{eq:bieab1}
\int_{-1}^1 G(s,t) \phi(t) dt = f(s) , \quad s \in (-1,1),
\end{align}
admits a solution $\phi \in \widetilde{H}^{-\frac{1}{2}}(-1,1)$, then the expansion $\phi$ as $\phi = \sum_{l=0}^\infty a_l w^{-1}T_l$ holds, with
$$
|a_l| \leq C \varrho^{-l},\quad \forall \ l \in \mathbb{N},
$$
for some $\varrho >1$.
\end{lemma}

\begin{proof}
If the solution exists, the expansion is possible since $\text{span}\{ w^{-1}T_n , n \in \mathbb{N}_0 \}$ is dense in $\widetilde{H}^{-\frac{1}{2}}(-1,1)$ (see \cite[Lemma C.2]{JP20}). Moreover, the norm of this space can be represented as
\begin{align}
\label{eq:norm1}
\| \phi\|_{\widetilde{H}^{-\frac{1}{2}}(-1,1)}^2 =\sum_{l=0}^\infty (1+l^2)^{-\frac{1}{2}} |{a}_l|^2 < \infty.
\end{align}
Since $R$ is analytic, it admits an expansion in terms of Chebyshev polynomials~{\cite[Theorem 8.1]{T13}} of the form:
\begin{align*}
R(s,t) = \sum_{p=0}^\infty \sum_{q=0}^\infty R_{p,q}T_p(s) T_q(t) \quad \text{with} \quad |R_{p,q}| \leq C \varrho^{-\max\{p,q\}}.
\end{align*}
On the other hand, by \cite[Lemma 4.14]{JP20}, it holds that
\begin{align}
\label{eq:bcofs}
\log|s-t|(J(s,t)-J_0) = \sum_{p=0}^\infty \sum_{q=0}^\infty B_{p,q} T_p(t) T_q(s) \quad \text{with} \quad |B_{p,q}| \leq C \min \{p^{-3},q^{-3}
\}.
\end{align}
Now, by combining the results  for the Laplace (see Remark~\ref{laplace}) case and the above expansions, we find that the integral equation \eqref{eq:bieab1} can be recasted as the following system of equations for the unknown coefficients $a_l$:
\be
\label{lsys}
J_0d_l a_l + \sum_{m=0}^\infty B_{l,m} a_m + \sum_{m=0}^\infty R_{l,m} a_m  = {f}_l, \quad l \in \mathbb{N}_0,
\en
where coefficients $d_l$ are those in \eqref{eq:lapbie} and ${f}_l$ denotes the $l$-th Chebyshev coefficient of ${f}$. Since ${f}$ is analytic, the coefficients ${f}_l$ decay exponentially fast. Also, it is clear that the third term on the left-hand side of (\ref{lsys}) decays exponentially. Thus, we have that
$$
|J_0 d_l a_l| - \left\vert  \sum_{m=0}^\infty B_{l,m} a_m  \right\vert \leq\left\vert J_0d_l a_l + \sum_{m=0}^\infty B_{l,m} a_m  \right\vert \leq C \varrho^{-l},
$$
and we conclude that either both terms $|j_0 d_l a_l|, \left\vert  \sum_{m=0}^\infty B_{l,m} a_m  \right\vert$ decay exponentially or they need to have the same order of decay (plus an exponentially decaying term). Assume first that they have the same decay order, by \eqref{eq:bcofs} then it holds that
\begin{align*}
\left\vert  \sum_{m=0}^\infty B_{l,m} a_m  \right\vert \leq l^{-3\beta} \sum_{m=0}^\infty
(m+1)^{-3\alpha} |a_m|
\end{align*}
for any real $\alpha, \beta >0$ such that $\alpha + \beta =1$. Furthermore, let $\alpha_1, \alpha_2>0$ such that $\alpha_1 + \alpha_2 = \alpha$. Thus, by the Cauchy-Schwarz inequality, one has
\begin{align*}
\left\vert  \sum_{m=0}^\infty B_{l,m} a_m  \right\vert^2 \leq l^{-6\beta} \sum_{m=0}^\infty
(m+1)^{-6\alpha_1} \sum_{m=0}^\infty (m+1)^{-6\alpha_2}|a_m|^2.
\end{align*}
If we select $\alpha_2 = \frac{1}{6}$, the right-most summation term can be bounded by $\|\phi\|_{\widetilde{H}^{-\frac{1}{2}}(-1,1)}$, the middle term $\sum_{m=0}^\infty
(m+1)^{-6\alpha_1}$ is finite if $\alpha_1 > \frac{1}{6}$. Therefore, by setting $\alpha=(1+\epsilon)/3$ and $\beta=(2-\epsilon)/3$ for some $0<\epsilon<2$, we get
\begin{align*}
\left\vert  \sum_{m=0}^\infty B_{l,m} a_m  \right\vert \leq C l^{-(2-\epsilon)} .
\end{align*}
Hence, since we assumed that both terms $|j_0 d_l a_l|, \left\vert  \sum_{m=0}^\infty B_{l,m} a_m  \right\vert$ have the same decay order, it can be concluded from the fact $d_l=\frac{\pi^2}{2l}, l\ge 1$ that $|a_l| \leq C l^{-1+\epsilon}$ which, however, further implies that
\begin{align*}
\left\vert  \sum_{m=0}^\infty B_{l,m} a_m  \right\vert \leq C l^{-(3-\epsilon)} .
\end{align*}
Consequently, both terms $|j_0 d_l a_l|, \left\vert  \sum_{m=0}^\infty B_{l,m} a_m  \right\vert$ can not decay with the same order, and so they both must decay exponentially in the $l$ variable.
\end{proof}
The following result generalizes the previous result for a more complicated form associated with the elastic hyper-singular BIE (\ref{NBIE}) for which the hyper-singular BIO can be reformulated as a combination of weakly-singular integrals and tangential derivatives.

\begin{lemma}
\label{lemma:coefsgeneralbies2}
For an analytic function $g$ defined over $(-1,1)$, we consider the hyper-singular integral equation:
\begin{eqnarray}
\label{eq:bieab2}
&&\frac{d}{ds}\int_{-1}^1 G^1(s,t) \psi^\prime(t) dt+
\frac{d}{ds}\int_{-1}^1 G^2(s,t) \psi(t) dt\\
&&+\int_{-1}^1 G^3(s,t) \psi^\prime(t) dt
 + \int_{-1}^1 G^4(s,t) \psi(t) dt = g(s)\nonumber
\end{eqnarray}
where the kernels $G^1, G^2,G^3,G^4$ can be expanded analogously to \eqref{eq:kersplit1}. If the integral equation (\ref{eq:bieab2}) admits a solution $\psi \in \widetilde{H}^{\frac{1}{2}}(-1,1)$, then it can be expanded as $\psi = \sum_{l=0}^\infty b_l wU_l$, and we also have that
$$
|b_l| \leq C \varrho^{-l},\quad \forall l \in \mathbb{N},
$$
for some $\varrho >1$.
\end{lemma}

\begin{proof}
We proceed as in the proof of Lemma \ref{lemma:coefsgeneralbies1}. Following to the kernel splitting (\ref{eq:kersplit1}), let $R^k, J^k, J_0^k$ denote the components of the decomposition of $G^k$ for $k = 1,2,3,4$, and denote by $R^k_{p,q}$, $B^k_{p,q}$, $p,q\in\N_0$ the coefficients of the Chebyshev polynomial expansions of $R^k$ and ${\log|s-t|(J^k(s,t)-J_0^k)}$ for $k=1,2,3,4$, respectively. The expansion of the solution $\psi$ in terms of the second-kind weighted Chebyshev polynomials follows by density, and the integral equation \eqref{eq:bieab2} is equivalent to the following system of equations for the coefficients $b_l$:
\ben
&& J^1_0 (l+1) d_{l+1} b_l+ \sum_{m=0}^\infty (m+1) B^1_{l+1,m+1}b_l + \sum_{m=0}^\infty (m+1) R^1_{l+1,m+1}b_l\\
&& +
\sum_{m^0}^\infty V^2_{l,m}b_m+
\frac{1}{2} \sum_{m=0}^\infty (R^2_{l+1,m}-R^2_{l+1,m+2})b_m+
\frac{1}{2} \sum_{m=0}^\infty (B^2_{l+1,m}-B^2_{l+1,m+2})b_m\\
&& + \sum_{m=0}^\infty V^3_{l,m}b_m+
\sum_{m=0}^\infty \frac{m+1}{2(l+1)} (R^3_{l+2,m+1}-R^3_{l,m+1})b_m\\
&&+\sum_{m=0}^\infty \frac{m+1}{2(l+1)} (B^3_{l+2,m+1}-B^3_{l,m+1})b_m \\
&& +
\sum_{m=0}^\infty V^4_{l,m}b_m +
\frac{1}{4(l+1)}\sum_{m=0}^\infty b_m (B^4_{l,m}-B^4_{l,m+2}-B^4_{l+2,m}+B^4_{l+2,m+2})\\
&& +\frac{1}{4(l+1)}\sum_{m=0}^\infty b_m (R^4_{l,m}-R^4_{l,m+2}-R^4_{l+2,m}+R^4_{l+2,m+2}) = \frac{1}{l+1}\widehat{g}_l, \quad l \in \mathbb{N}_0,
\enn
where $\widehat{g}_l$ corresponds to the $l$th Chebyshev coefficient of the second kind (i.e.~$\widehat{g}_l := \int_{-1}^1 g(t) w(t)U_l(t) dt$), and
\ben
&& V^2_{l,m} =  \frac{1}{2}J_0^2d_{l+1}(\delta_{l-1,m} -\delta_{l+1,m})\le Cl^{-1},\\
&& V^3_{l,m} = \frac{1}{2(l+1)}J_0^3(d_{l+2}(l+2)\delta_{l+1,m}-d_l l\delta_{l-1,m})\le Cl^{-1},\\
&& V^4_{l,m} = \frac{1}{4(l+1)}J_0^4(2d_l\delta_{l,m}- d_{l-2}\delta_{l-2,m}-d_{l+2}\delta_{l+2,m})\le Cl^{-2}.
\enn
Since $g$ is assumed to be analytic, the coefficients $\widehat{g}_l$ decay exponentially. Note that all the summation terms involving $R^k_{p,q}$, $k=1,2,3,4$, $p,q\in\N_0$ decay exponentially. Thus, the first term $J^1_0 (l+1) d_{l+1} b_l$ plus the summation terms involving coefficients $B^k_{p,q}$, $V^k_{p,q}$, $k=1,2,3,4$, $p,q\in\N_0$ should also decay exponentially. Similarly to the argument presented in Lemma~\ref{lemma:coefsgeneralbies1} for the case of weakly-singular integral equation, it can be concluded that the the summation of $J_0^1(l+1) d_{l+1}b_l$ plus the terms involving the coefficients $V^k_{p,q},\ k=2,3,4,\ p,q\in\N_0$ must also decay exponentially. Therefore, the decay properties of $d_l$ and $V^k_{l,m},\ k=2,3,4$ implies that the coefficients $b_l$ have to decay exponentially in $l$.
\end{proof}

\subsection{Convergence results}
\label{sec:3.3}

In order to use the abstract results from Lemmas~\ref{lemma:coefsgeneralbies1}--\ref{lemma:coefsgeneralbies2} for the BIEs appearing in Problem \ref{problem1}, we first need to verify that they can be recasted as integral equations on $(-1,1)$, and also that the corresponding kernels can be decomposed as in \eqref{eq:kersplit1}, and that the resulting right-hand sides are analytic.

To start with, notice that the general structure of the weakly-singular BIE in Problem \ref{problem1} is
\be
\label{wsbie}
\sum_{j=1}^M\int_{\Gamma_j} G(\bs x,\bs y) \varphi_j(\bs y) d \bs y  =
f_{i}(\bs x),  \quad \bs x \in \Gamma_i,\ i = 1,\hdots, M,
\en
where $G(\cdot,\cdot)$ is a generic kernel corresponds to the four components of the fundamental solution $\mathbb{E}(\cdot,\cdot)$ which can be expressed as
\ben
\label{eq:fsolV}
{\bb E}(\bs x,\bs y)&=&  \frac{i}{4\mu}H_0^{(1)}(\kappa_s|\bs x-\bs y|){\bb I}\nonumber\\
&-& \frac{i}{4\rho\omega^2|\bs x-\bs y|} \left[\kappa_sH_1^{(1)}(\kappa_s|\bs x-\bs y|)-\kappa_pH_1^{(1)}(\kappa_p|\bs x-\bs y|)\right]{\bb I} \nonumber\\
&+& \frac{i(\bs x-\bs y)(\bs x-\bs y)^\top}{4\rho\omega^2|\bs x-\bs y|^2} \left[\kappa_s^2H_2^{(1)}(\kappa_s|\bs x-\bs y|)-\kappa_p^2H_2^{(1)}(\kappa_p|\bs x-\bs y|)\right]
\enn
where $H^{(1)}_1(\cdot)$ and $H^{(1)}_2(\cdot)$ are first-kind Hankel functions of first and second order. Employing the parametrizations $\bm {r}_i$ of each open-arc $\Gamma_i$ for $i=1, \hdots,M$, the integral equation (\ref{wsbie}) can be rewritten as
\be
\label{eq:bieabs3}
\sum_{j=1}^M\int_{-1}^1 G(\bm {r}_i(s),\bm{r}_j(t)) \widehat{\varphi}_j(t) dt  = f_i\circ \mathbf{r}_i(s),  \quad s \in (-1,1),\ i = 1,\hdots, M ,
\en
where $\widehat{\varphi}_j = \varphi_j \circ \bm {r}_j \|  \bm{r}_j' \|$.
Hence, under the assumption that the right-hand sides of Problem \ref{problem1} and the parametrizations are analytic, the right-hand side of \eqref{eq:bieabs3} is also analytic and the hypotheses of the previous lemma as to the right-hand side are fulfilled. On the other hand, for $i \neq j$, the components of ${\bb E}(\bm r_i(s),\bm r_j(t))$ are analytic functions for which the corresponding terms in the decomposition (\ref{eq:kersplit1}) would be  $J_0 = J(s,t) = 0$, while for $i= j $, the decomposition of the form \eqref{eq:kersplit1} can be obtained from the series expansion of Bessel functions (see \cite[Chapter 9]{stg}).

Next, we consider the hyper-singular BIE in Problem \ref{problem1}. It follows from the regularization technique presented in \cite{YHX17} that the hyper-singular BIO $\mathcal{W}_{ij}$ admits the equivalent form:
\ben
\label{eq:fsolW}
&&\mathcal{W}_{ij}[{\bs \psi}_j](\bs x)\nonumber=\\
&& -\rho\omega^2 \int_{\Gamma_j}\left[\gamma_{k_s}({\bs x},{\bs y})(2\bs\nu_{\bs x}\bs\nu_{\bs y}^\top- \bs\nu_{\bs y}\bs\nu_{\bs x}^\top-\bs\nu_{\bs x}^\top\bs\nu_{\bs y}\mathbb{I}) - \gamma_{k_p}({\bs x},{\bs y})\bs\nu_{\bs x}\bs\nu_{\bs y}^\top\right]{\bs \psi}_j({\bs y})ds_{\bs y}\nonumber\\
&&+ 4\mu^2\dfrac{d}{ds_{\bs x}}\int_{\Gamma_j}
\left[\mathbb{A}\mathbb{E}({\bs x},{\bs y})\mathbb{A}+ \frac{1}{\mu}\gamma
_{k_s}({\bs x},{\bs y})\mathbb{I} \right]\dfrac{d{\bs \psi}_j({\bs y})}{ds_{\bs y}} ds_{\bs y}\nonumber\\
&&- 2\mu \int_{\Gamma_j}  \bs\nu_{\bs x} \nabla _{\bs x}^\top [\gamma_{k_s}({\bs x},{\bs y})-\gamma_{k_p}({\bs x},{\bs y})]\mathbb{A}\dfrac{d{\bs \psi}_j({\bs y})}{ds_{\bs y}}ds_{\bs y}\nonumber\\
&&- 2\mu\frac{d}{ds_{\bs x}}\int_{\Gamma_j} \mathbb{A}\nabla_{\bs y}
[\gamma_{k_s}({\bs x},{\bs y})-\gamma_{k_p}({\bs x},{\bs y})]\bs\nu_{\bs y}^\top {\bs \psi}_j({\bs y})ds_{\bs y},
\enn
with $$\quad \mathbb{A}=\begin{bmatrix}
0 & -1\\
1 & 0
\end{bmatrix}.$$
Then the analyticity assumptions of both right-hand sides of Problem \ref{problem1}, the parametrizations of open-arcs together with the series expansion of Bessel functions~\cite[Chapter 9]{stg} imply that the hyper-singular BIE in Problem \ref{problem1} can also be expressed as integral equations on $(-1,1)$ and the conditions of Lemma~\ref{lemma:coefsgeneralbies2} are satisfied.

\begin{corollary}
\label{cor:elasticitycofs}
Let $\bs \phi\in\widetilde{\mathbb{H}}^{-\half}(\Gamma)$ and $\bs \psi\in\widetilde{\mathbb{H}}^{\half}(\Gamma)$ be the unique solutions of Problem \ref{problem1}, then they admit the decompositions:
\ben
&&{\bs \phi}_i \circ \mathbf{r}_i(t) = \sum_{p=1}^2 \sum_{l=0}^\infty (a_i)_{l,p} w^{-1}(t)T_l(t){\bs e}_p , \quad i=1,\hdots,M, \\
&&{\bs \psi}_i \circ \mathbf{r}_i(t) = \sum_{p=1}^2 \sum_{l=0}^\infty (b_i)_{l,p} w(t)U_l(t){\bs e}_p, \quad i=1,\hdots,M,
\enn
and the following bounds on the coefficients hold
\begin{align*}
| (a_i)_{l,p}| \leq C \varrho^{-l}, \quad
| (b_i)_{l,p}| \leq C \varrho^{-l},\quad l\in\N,\ p=1,2,\ i=1,\hdots,M.
\end{align*}
\end{corollary}
\begin{proof}
We will only prove for $\bs \phi$ as for $\bs \psi$ the arguments are similar but using Lemma \ref{lemma:coefsgeneralbies2} instead of Lemma \ref{lemma:coefsgeneralbies1}.

Since $\bs \phi \in \widetilde{\mathbb{H}}^{-\half}(\Gamma)$, each component $\bs \phi_i$ is in $\widetilde{H}^{-\half}(\Gamma_i)^2$. Furthermore, one can directly show---using duality arguments and the Sobolev spaces definition via the Sobolev-Slobodeckii norm (see \cite[Chapter 2]{M00})---that for every $v \in \widetilde{H}^{-\half}(\Gamma_i)$, it holds that
$$
\| v \|_{\widetilde{H}^{-\half}(\Gamma_i)} \cong
 \| v \circ \bs r_i \|_{\widetilde{H}^{-\half}(-1,1)}.
$$
Hence, $\bs \phi_i \circ \bs r_i \in \widetilde{H}^{-\half}(-1,1)^2$, and therefore the expansion of each component follows from the density of weighed first-kind Chebyshev polynomials in $\widetilde{H}^{-\half}(-1,1)$. The corresponding BIE for $\bs \phi$ reads
\begin{align*}
\int_{-1}^1 \bb E(\bm r_i(s),\bm r_i(t)) \bs \phi_i \circ \bm r_i(t) dt + \sum_{i \neq j} \int_{-1}^1 \bb E(\bm r_i(s),\bm r_j(t)) \bs \phi_j \circ \bm r_j(t) dt = \bs f \circ \bm r_i(s),
\end{align*}
Since $\bb E(\bm r_i(s),\bm r_j(t))$ is analytic for $i \neq j$, when the above equations are transformed into a system of equation for the coefficients $(a_i)_{l,p}$ these terms do not alter the convergence rate. Thus, the result follows directly by Lemma \ref{lemma:coefsgeneralbies1}.
\end{proof}

With the exponentially decaying properties of the sequences $(a_i)_{l,p}$, $(b_i)_{l,p}$, $l\in\N$ for $i=1,\hdots,M$, $p=1,2$, one can easily obtain the convergence rate of our Galerkin method presented in Section~\ref{sec:3.1}.

\begin{corollary}
\label{cor:convergencerate}
Let $\bs \phi\in\widetilde{\mathbb{H}}^{-\half}(\Gamma)$ and $\bs \psi\in\widetilde{\mathbb{H}}^{\half}(\Gamma)$ be the unique solutions of Problem \ref{problem1}. Let $N, N_0 \in \mathbb{N}$, with $N>N_0$ be such that Problem \ref{problem2} has unique solutions, and $\bs \phi^N, \bs \psi^N$ denote the corresponding discrete approximations (defined as in \eqref{eq:discretesols}). Then, there exists $\varrho>1$ and a constant $C >0$ such that
\begin{align*}
\|\bs \phi-\bs \phi^N\|_{\widetilde{\mathbb{H}}^{-\half}(\Gamma)} &\leq C \varrho^{-N},\\
\|\bs \psi-\bs \psi^N\|_{\widetilde{\mathbb{H}}^{\half}(\Gamma)} &\leq C \varrho^{-N}.
\end{align*}
\end{corollary}
\begin{proof} Since we are using a Galerkin discretization of a coercive problem, such value of $N$ exists (cf.~\cite[Section 4.2]{SS10} and Lemma \ref{lemma:qopty}) possibly differing for Dirichlet and Neumann cases. Once again we focus only on the Dirichlet case as the Neumann one follows verbatim. From the quasi-optimality result (\ref{quasi1})-(\ref{quasi2}) and the norm equivalence used in the proof of Corollary \ref{cor:elasticitycofs}, we have that
\ben
\|\bs\phi-\bs\phi^N\|_{\widetilde{\mathbb{H}}^{-\half}(\Gamma)} \leq C \sum_{i=1}^M \sum_{p=1}^2 \inf_{\varphi_N \in \widehat{\mathbf{T}}_N}\left\| \sum_{l=0}^\infty (a_i)_{l,p} w^{-1}(t)T_l(t)-\varphi_N \right\|_{\widetilde{H}^{-\frac{1}{2}}(-1,1)},
\enn
where $\widehat{\mathbf{T}}_N := \text{span}\{ w^{-1}T_n, n=0,\hdots,N\}$. From Corollary \ref{cor:elasticitycofs}, for fixed $i \in \{1,..,M\}$, and $p \in \{1,2\}$, we can choose $\varphi_N =  \sum_{l = 0}^N (a_i)_{l,p} w^{-1}T_l$. Thus, we derive the following bound
\begin{align*}
\|\bs \phi-\bs \phi^N\|_{\widetilde{\mathbb{H}}^{-\half}(\Gamma)} \leq C \sum_{i=1}^M \sum_{p=1}^2 \left\Vert\sum_{l >N} (a_i)_{l,p} w^{-1}T_l \right\Vert_{\widetilde{H}^{-\frac{1}{2}}(-1,1)}.
\end{align*}
The right-most term can be bounded as follows
\begin{align*}
\left\| \sum_{l >N} (a_i)_{l,p}w^{-1}T_l \right\|^2_{\widetilde{H}^{-\frac{1}{2}}(-1,1)} = \sum_{l>N} (1+l^2)^{-\frac{1}{2}} |(a_i)_{l,p}|^2 \leq C \sum_{l>N} (1+l^2)^{-\frac{1}{2}} \varrho^{-2l},
 \end{align*}
where the last inequality follows from the bounds in Corollary \ref{cor:elasticitycofs}. The final result then follows directly by recalling the formula for geometric sums.
\end{proof}

\begin{remark}
Though the singular edge behavior was explicitly included in the discrete spaces of Section \ref{sec:3.1}, obtaining convergence rates does not require particular assumptions on solutions $\bs \phi$, $\bs \psi$ singularities. In fact, we can obtain as a corollary---arguing as in \cite[Theorem 8.3]{T13}---that the components of the solutions of the Dirichlet problem, mapped back to $[-1,1]$, could be written as $h(t)(1-t^2)^{-\frac{1}{2}}$, where $h:[-1,1]\rightarrow \mathbb{R}$, is an analytic function. Similarly, solutions of the Neumann problem have the general form $h(t) (1-t^2)^{\frac{1}{2}}$, with $h$ as before.
\end{remark}

\section{Numerical implementation and experiments}
\label{sec:4}

Before performing numerical experiments, we provide implementation details concerning the fast computation of matrix entries improving also the strategies presented in {\cite[Section 6]{JP20}}.

\subsection{Implementation strategy}
\label{sec:4.1}

Following the definition of discrete spaces in Section \ref{sec:3.1} and the fundamental solution representation in Section \ref{sec:3.3}, the numerical implementation of the method relies on computing integrals of the generic forms:
\begin{align*}
I^1_{l,m}&:=\int_{-1}^1 \int_{-1}^1 G(\bm r_i(s), \bm r_j(t)) \varphi_m(s) \varphi_l(t) ds dt, \\
I^2_{l} &:= \int_{-1}^1 h(t) \varphi_l(t)dt,
\end{align*}
where $l,m \in \{0,\hdots,N\}$, and $N$ is the parameter determining the number bases per arc and proportional to the dimension of the discretization space. Galerkin formulations for the weakly- and hyper-singular BIOs imply that the functions $\varphi_l$ could take the one of the following structure in terms of Chebyshev polynomials:
$$
w^{-1}(t)T_l(t), \ w(t) U_l(t), \ \frac{d}{dt} w(t)U_l(t).
$$
Hence, from the identities
$$
w(t) U_l(t) = \frac{1}{2}w^{-1}(t) (T_l(t)-T_{l+2}(t)), \quad
\frac{d}{dt} w(t)U_l(t) = -(l+1) w^{-1}(t)T_{l+1}(t),
$$
we could reduce any computation to the case $\varphi_l (t) = w^{-1}(t)T_l(t)$.

For $I^1$, following the expansion of Hankel functions~\cite[9.1.13]{stg}, the term $G(\bm r_i(s), \bm r_j (t) )$ could be expressed as functions of the form:
$$
G(\bm r_i(s), \bm r_j (t) ) = \log |s-t| J(s,t) + R(s,t),
$$
where $J(\cdot, \cdot)$, $R(\cdot, \cdot)$ are analytic functions, and in particular $J(\cdot, \cdot) = 0$ if $i\neq j$. For $I^2$, the function $h(\cdot)$ is assumed to be analytic.

The implementation of the proposed spectral Galerkin method is then achieved using the techniques presented in \cite{JP20}. For $I^2$ we can find an interpolation approximation using the fast Fourier transform (FFT) of $h(\cdot)$:
$$h(t) \approx \sum_{n=0}^{N_c} h_n T_n(t),$$
with an error decaying exponentially for increasing values of $N_c$. Then, the integrals {$I^2$} are found using the orthogonality relations of Section \ref{sec:3.1}. In this work, we improve the implementation in \cite{JP20} by noticing that we can select $N_c$ adaptively. In fact, from \cite[Chapters 5 and 8]{T13}, we know that
$$
|h_n| \leq C \rho^{-n},
$$
for some $\rho>1$. The selection of $N_c$ can be done in two stages:
\begin{enumerate}
\item
Starting from $N_c = 1$ we compute the sequence $\{h_n\}_{n=0}^{N_c}$, and check if the last two coefficients are smaller than a given tolerance (typically {\texttt{tol}}=$10^{-12}$), if not we doubles the value of $N_c$.
\item
With a value of $N_c$ that ensure that the the last entries of $h_n$ are smaller than {\texttt{tol}}, we use a bisection search between $\frac{N_c}{2}$ and $N_c$ for the minimum value that still give as that the last two entries are smaller than the given tolerance.
\end{enumerate}
Furthermore, from the orthogonality properties of the Chebyshev polynomials, $I^2_l$ would be proportional to $h_l$. It is only necessary to compute the entries of $I^2_l$ that are bigger than the given tolerance, and thus we reduce memory requirements.

For $I^1$ the extension of the previous idea is direct. First, consider the case $i \neq j$, thus we only need to find the approximation
$$
R(t,s) \approx \sum_{p =0}^{N_c} \sum_{q =0}^{N_c} R_{p,q} T_p(t) T_q(s),
$$
which is done as in the previous case with the only difference that for finding $N_c$ we do not find the full bi-variate sequence $R_{p,q}$ but instead we use a greedy algorithm that selects $N_c$ from the following approximation
$$
R(t,0) \approx \sum_{p=0}^{N_c} R_p T_p(t).
$$
As in the case of $I^2$ this implementation would give us a sparse representation of the matrices as not all the entries are computed but only the one that the greedy algorithm estimates as bigger than the given tolerance.
Finally, we consider $I^1$ for $i=j$. The regular part $R(\cdot, \cdot)$ is integrated as in the case where $i \neq j$, thus we are left with the approximation of integrals of the form
$$
I^S_{l,m} = \int_{-1}^1 \int_{-1}^1 \log|t-s| J(t,s) w^{-1} T_m(s) w^{-1} T_l(t) ds dt.
$$
From \eqref{eq:logcofs}, we have that
$$
\log|t-s|  = \sum_{n=0}^\infty d_n T_n(t) T_n(s),
$$
for a known sequence $d_n$. Following the computations for the regular part we can construct the approximation
$$
J(t,s) \approx \sum_{p =0}^{N_c} \sum_{q =0}^{N_c} J_{p,q} T_p(t) T_q(s),
$$
Combining the last two equation and using the identity $T_a(t)T_b(t) = \frac{1}{2} (T_{a+b}(t)+T_{|a-b|}(t))$, we obtain
$$
I^S_{l,m} = \sum_{n=0}^\infty
\frac{d_n}{4} (J_{|l-n|,|m-n|}+J_{n+l,|m-n|}+J_{|l-n|,n+m}+J_{n+l,m+n}).
$$
The last sum is implicitly truncated as we assumed that $J_{p,q} =0$ if $p >N_c$ or $q> N_c$, hence the maximum number of terms in the sum is $N_c +N$. We remark that the computation of $I^S_{l,m}$ could be accelerated using convolution identities for discrete transform.
\subsection{Numerical results}
\label{sec:4.2}
We now present some numerical examples to illustrate our claims. Throughout unless is stated otherwise, we fix the parameter as $\mu = 1, \lambda = 2, \rho =1, \omega = 50$, and consider the scattering problems of a $p-$plane incident wave given by
$$
\bs P(\bs x) = \bs d e^{i k_p \bs x \cdot \bs d}
$$
where $\bs d = (\cos \alpha , \sin \alpha)^\top$, $\alpha\in[0,2\pi)$ being the incidence angle, and $k_p^2 =\dfrac{\omega^2 \rho}{\lambda+2\mu}$, as before. Then, the right-hand sides for the Dirichlet and Neumann problems are given by $\bs f_i=-\gamma_{D,i}\bs P$ and $\bs g_i=-\gamma_{N,i}\bs P$, respectively. Numerical errors shown as follows are defined by
\ben
\|\bphi_{\mathrm{num}}-\bphi_{\mathrm{ref}}\|_{\widetilde{\mathbb{H}}^{-\half}(\Gamma)} \quad\text{and}\quad
\|\bpsi_{\mathrm{num}}-\bpsi_{\mathrm{ref}}\|_{\widetilde{\mathbb{H}}^{\half}(\Gamma)}
\enn
for the Dirichlet and Neumann problems, respectively, wherein $\bphi_{\mathrm{num}}, \bpsi_{\mathrm{num}}$ are the numerical solutions and the reference solutions $\bphi_{\mathrm{ref}}, \bpsi_{\mathrm{ref}}$ are obtained as a numerical solution for sufficiently fine discretizations or overkill solutions.

\begin{figure}[htbp]
\centering
\begin{tabular}{cc}
\includegraphics[scale=0.07]{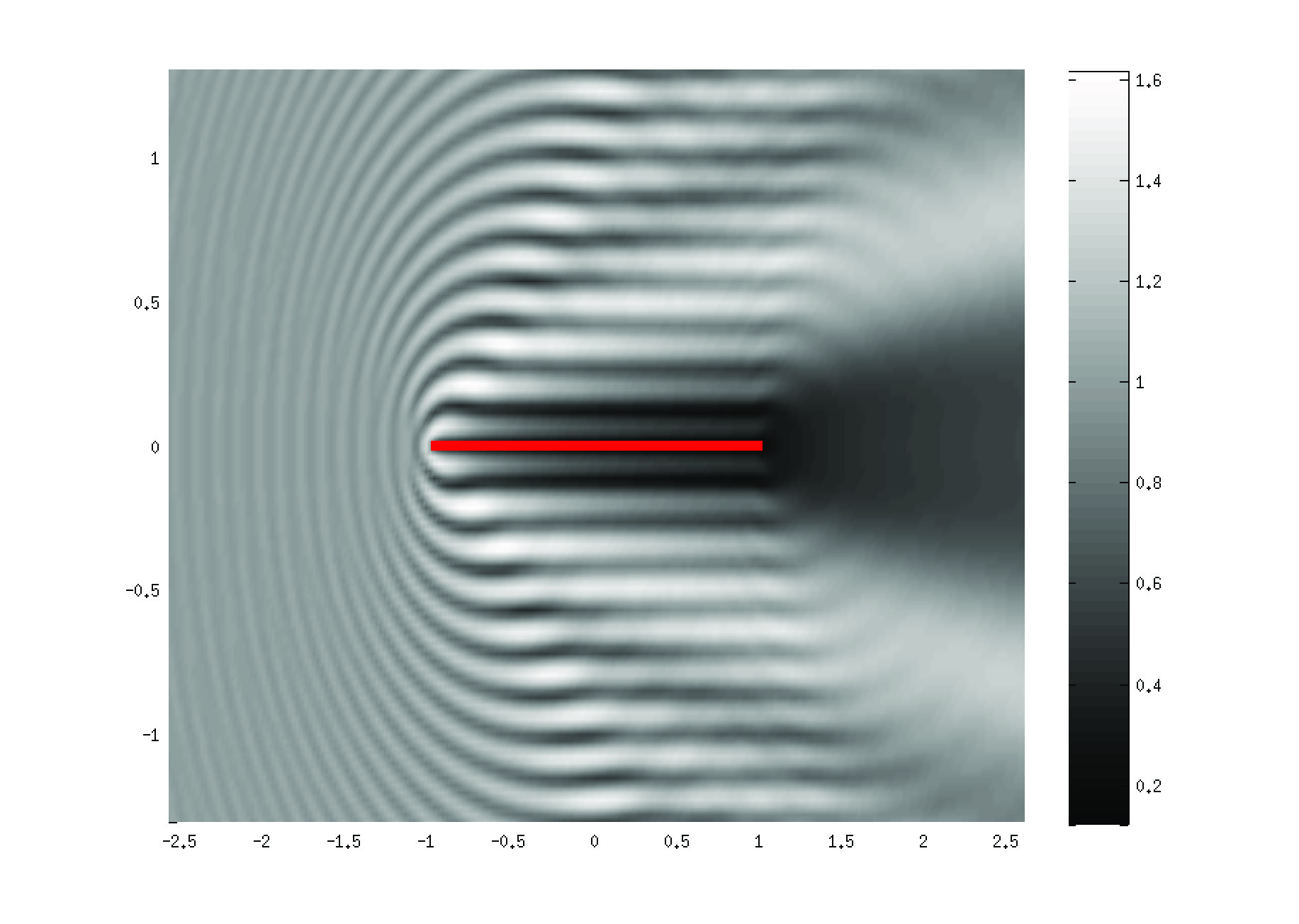} &
\includegraphics[scale=0.067]{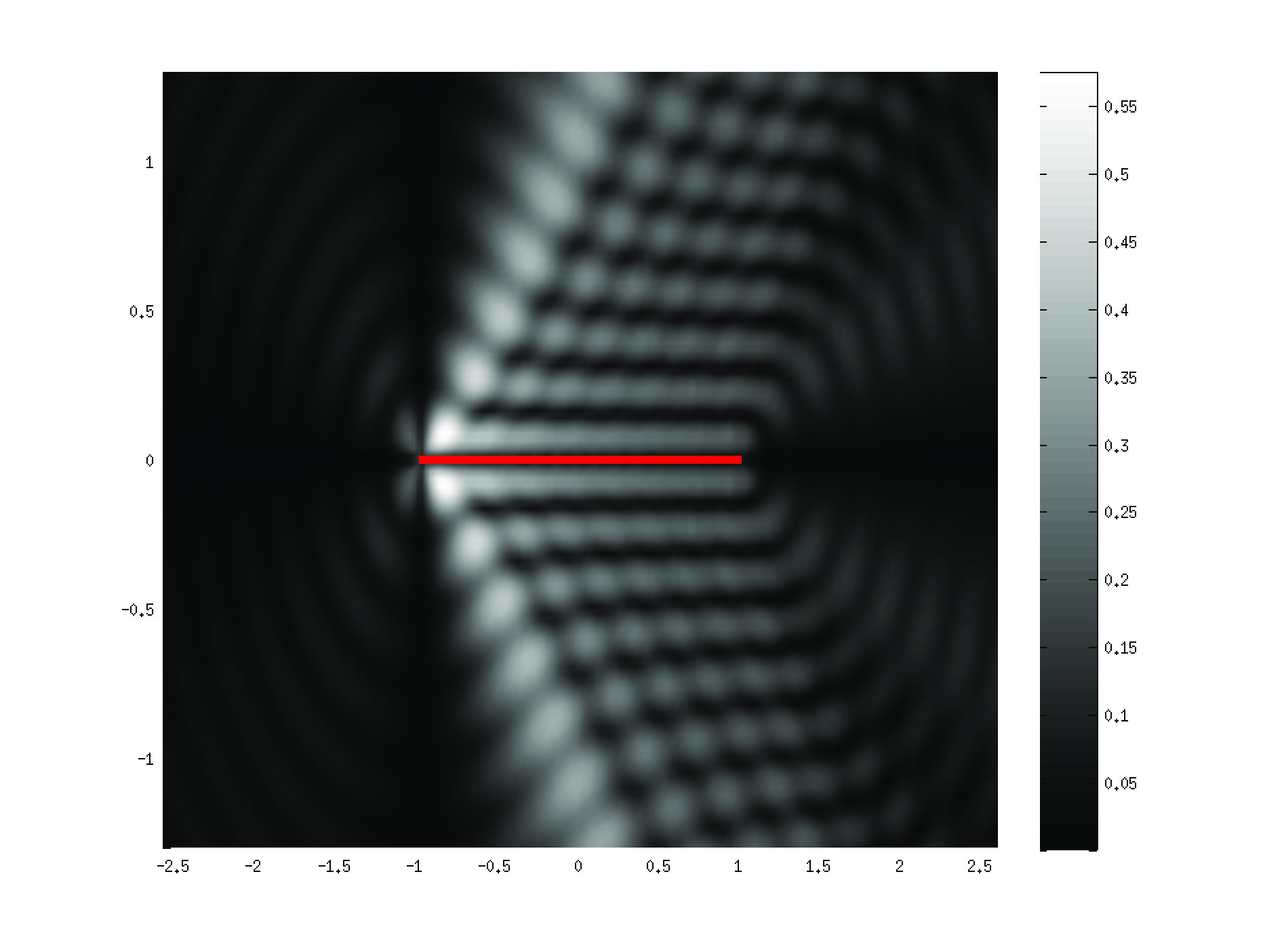} \\
(a) $U_1$ for Dirichlet Case & (b) $U_2$ for Dirichlet Case \\
\includegraphics[scale=0.29]{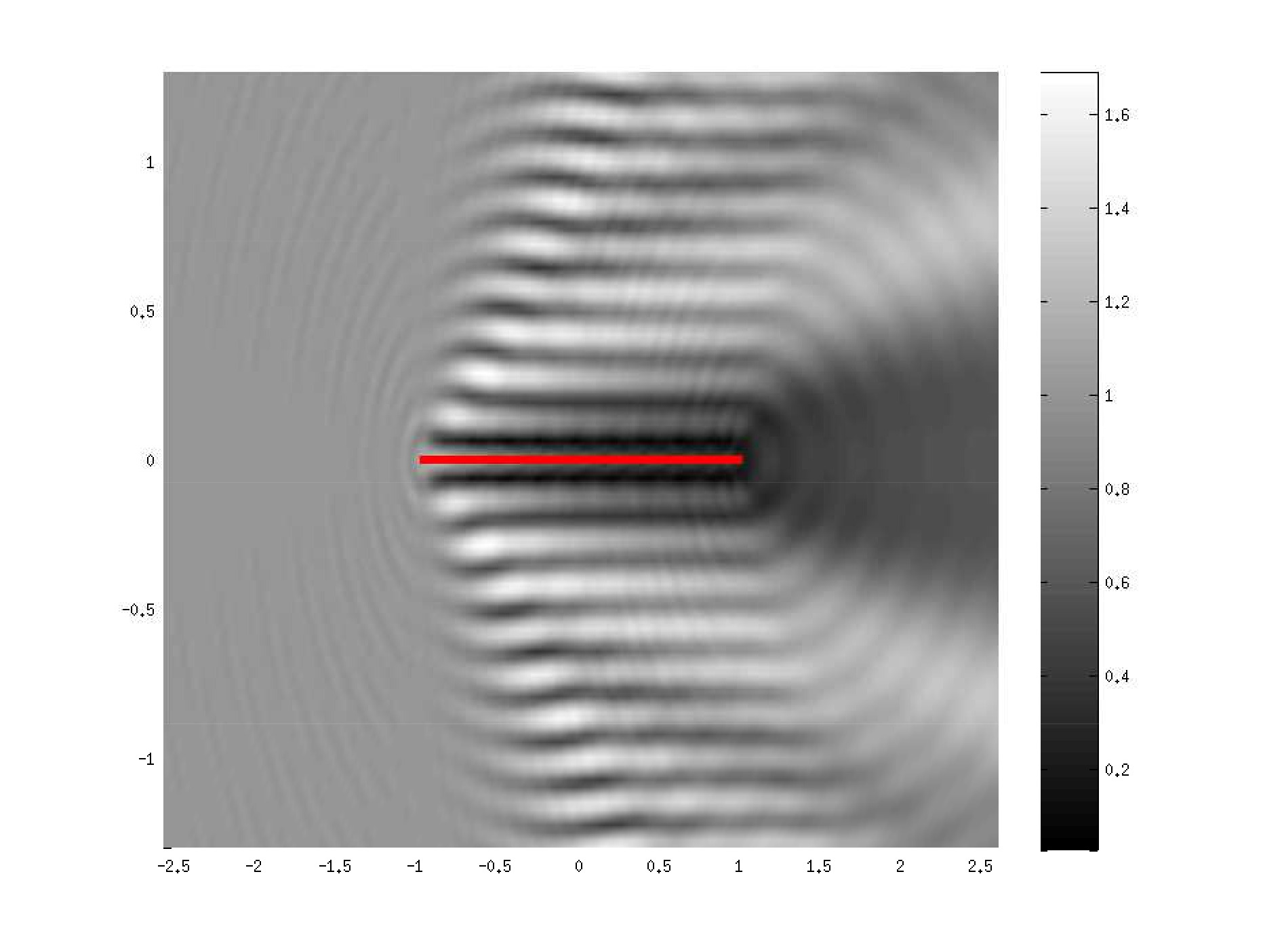} &
\includegraphics[scale=0.07]{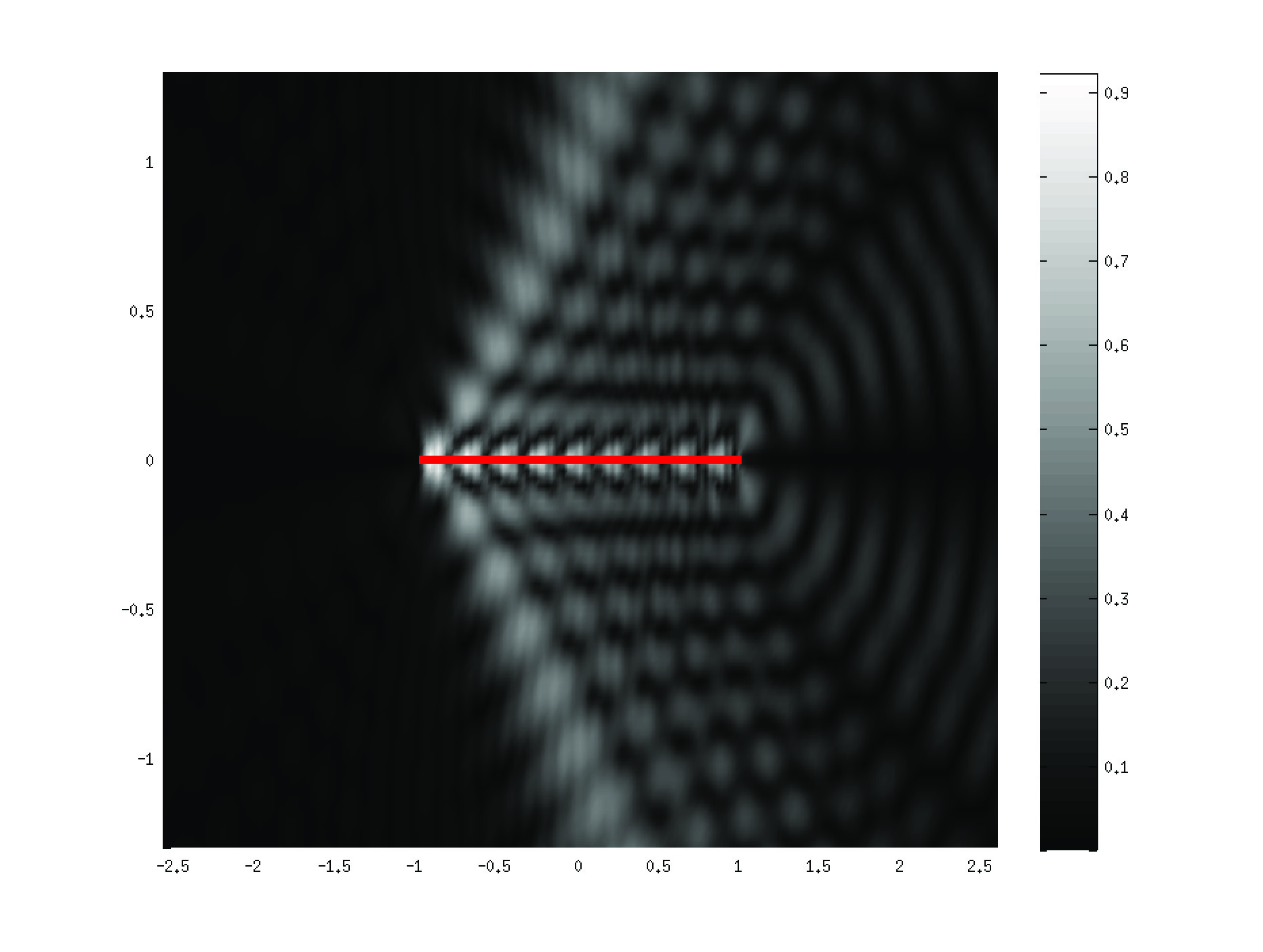} \\
(c) $U_1$ for Neumann Case & (d) $U_2$ for Neumann Case
\end{tabular}
\caption{Total field absolute values for the scattering problems on a single line segment (Example 1).}
\label{fig:linefields}
\end{figure}

\begin{figure}[htbp]
\centering
\includegraphics[scale=0.4]{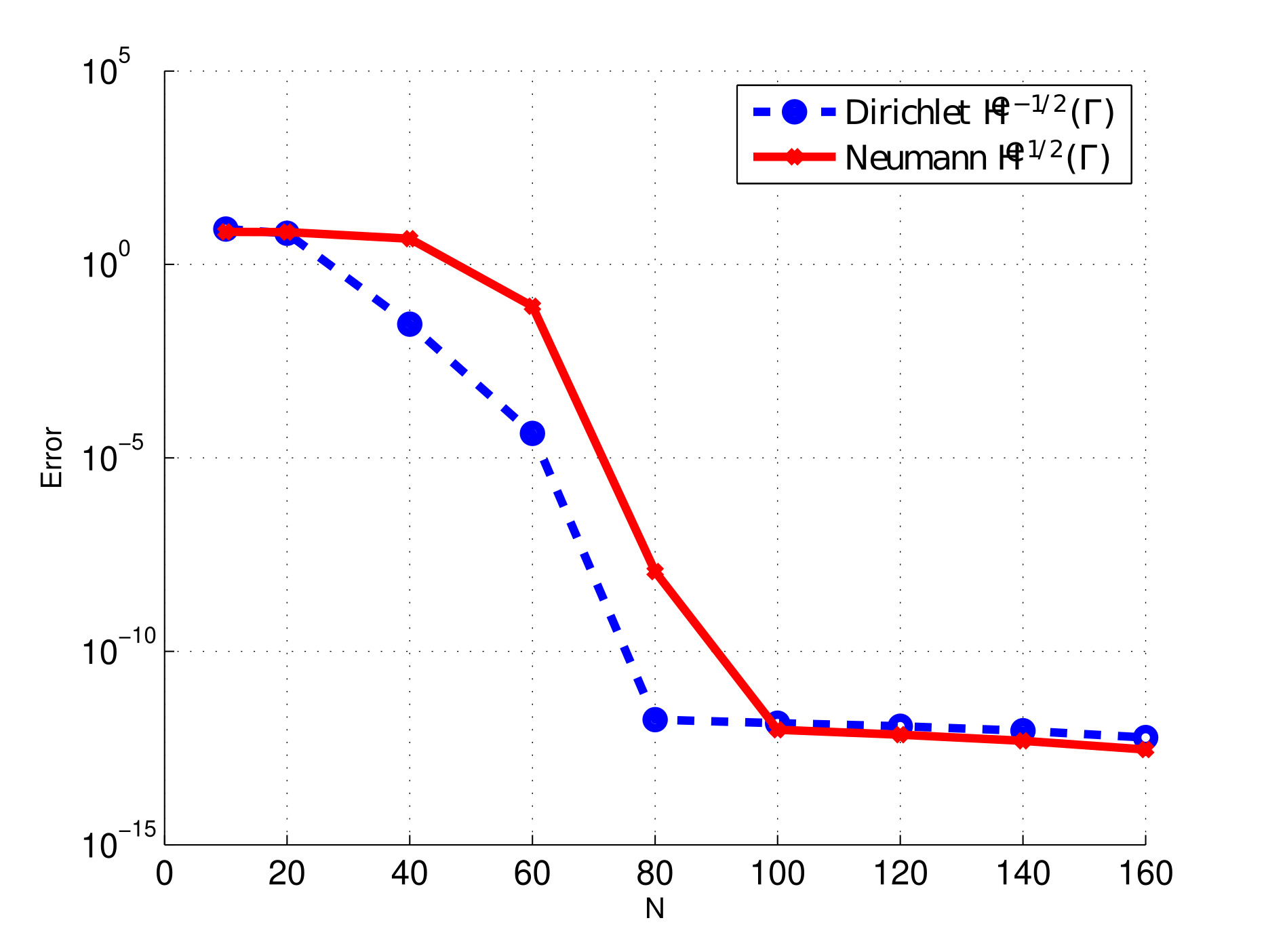}
\caption{Error convergence for the scattering problems on a single line segment (Example 1).}
\label{fig:errorline}
\end{figure}

{\bf Example 1.} We first consider the simple single line segment case, i.e., $\Gamma=\{\bs x=(t,0)\in\R^2: t\in[-1,1]\}$ and choose $\alpha = 0$. The total volume fields $\bs U^{\mathrm{tot}}=\bs U + \bs P$ for the Dirichlet and Neumann problems are displayed in Figure~\ref{fig:linefields}. Figure~\ref{fig:errorline} shows the numerical errors for an increasing number of polynomials basis for both problems, which demonstrates the exponential convergence of the proposed spectral Galerkin method in this setting.

\begin{figure}[htbp]
\centering
\begin{tabular}{cc}
\includegraphics[scale=0.07]{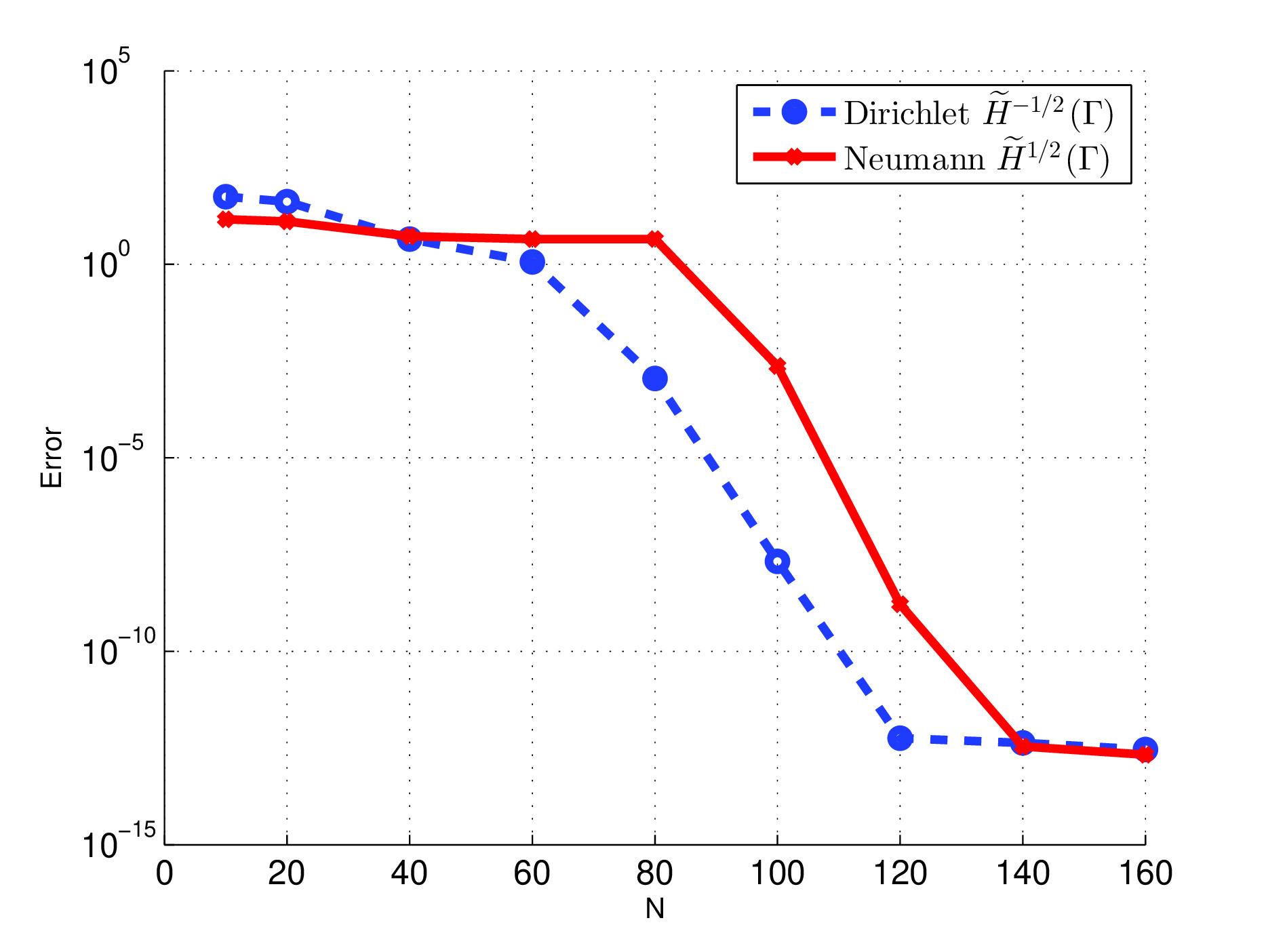} &
\includegraphics[scale=0.07]{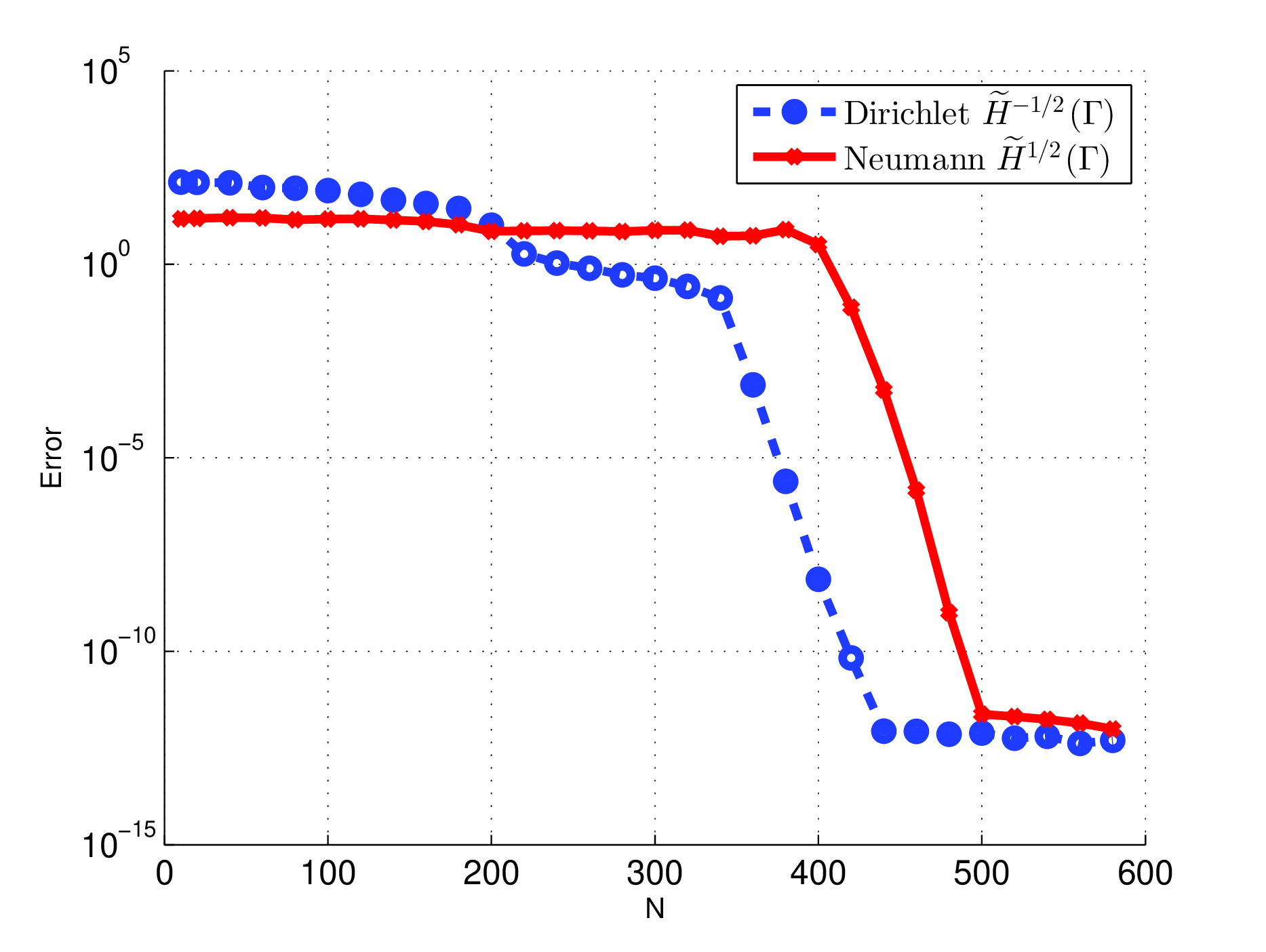} \\
(a)  &  (b)
\end{tabular}
\caption{Numerical errors for the scattering problems by a semi-circle (a) or a spiral-shaped (b) arc (Example 2).}
\label{fig:errorscircle}
\end{figure}

\begin{figure}[htbp]
\centering
\begin{tabular}{cc}
\includegraphics[scale=0.07]{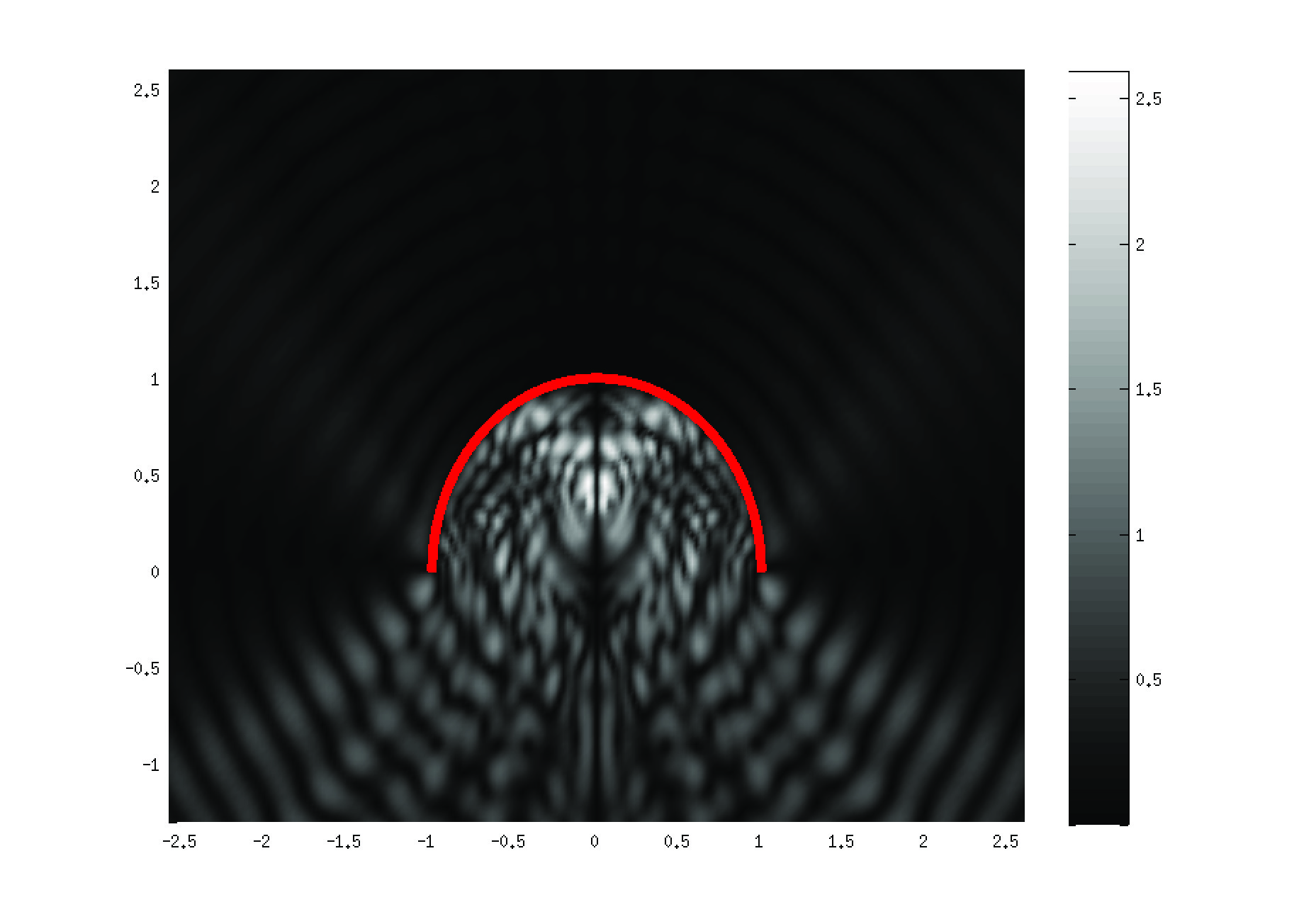} &
\includegraphics[scale=0.07]{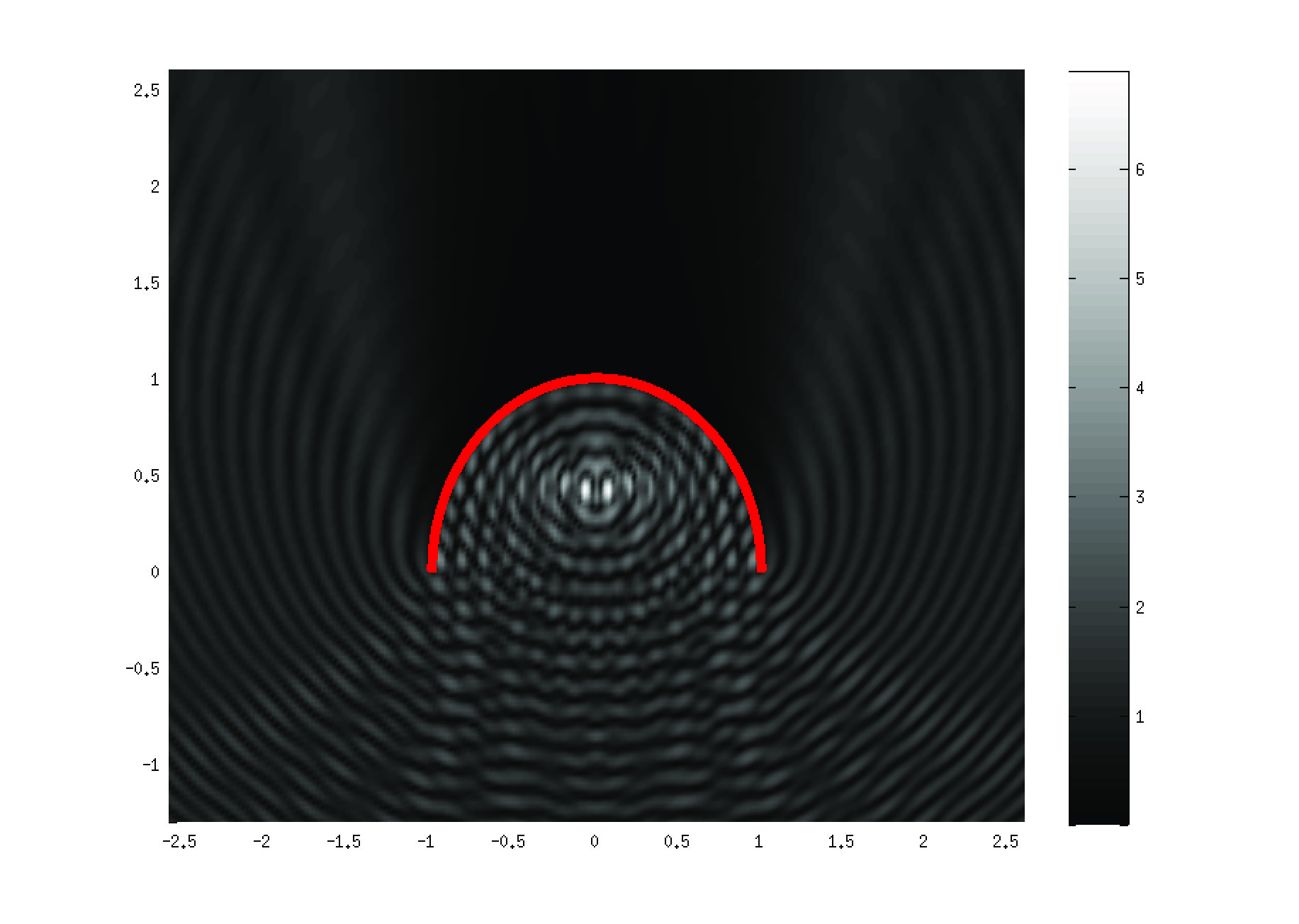} \\
(a) $U_1$ for Dirichlet Case & (b) $U_2$ for Dirichlet Case \\
\includegraphics[scale=0.07]{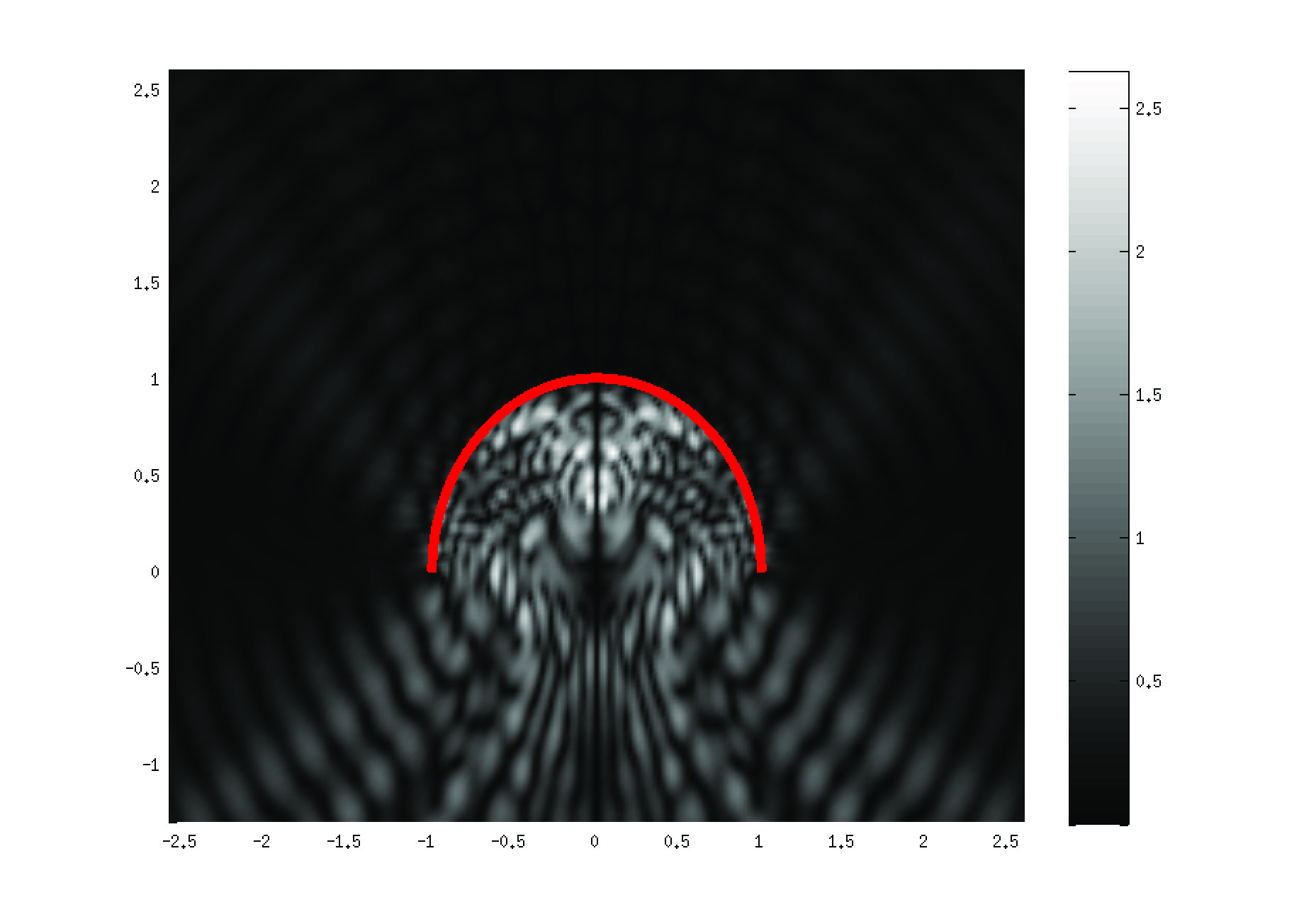} &
\includegraphics[scale=0.07]{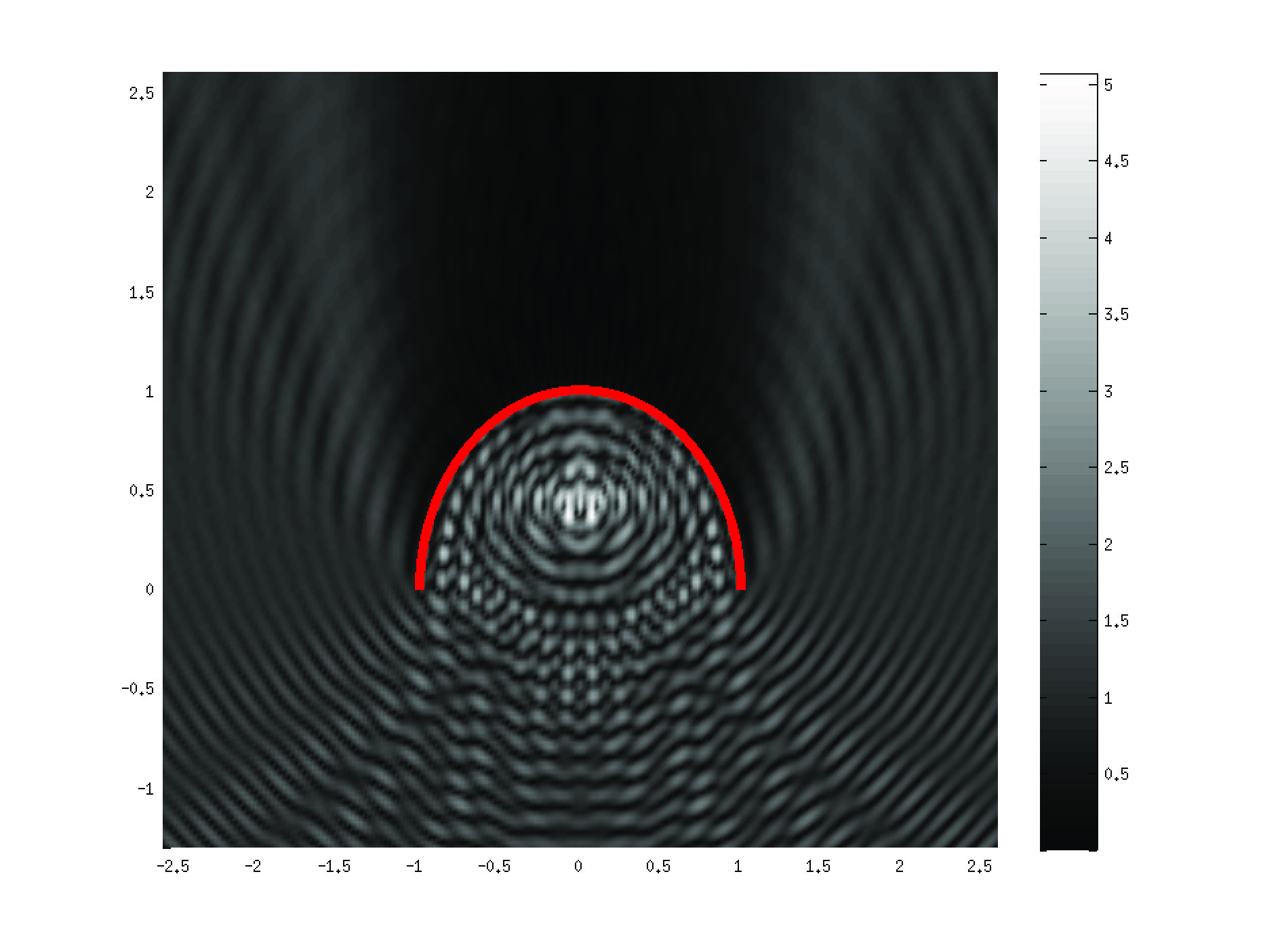} \\
(c) $U_1$ for Neumann Case & (d) $U_2$ for Neumann Case
\end{tabular}
\caption{Total field absolute values for the scattering problems by a semi-circle arc (Example 2).}
\label{fig:scirclefields}
\end{figure}

\begin{figure}[htbp]
\centering
\begin{tabular}{cc}
\includegraphics[scale=0.07]{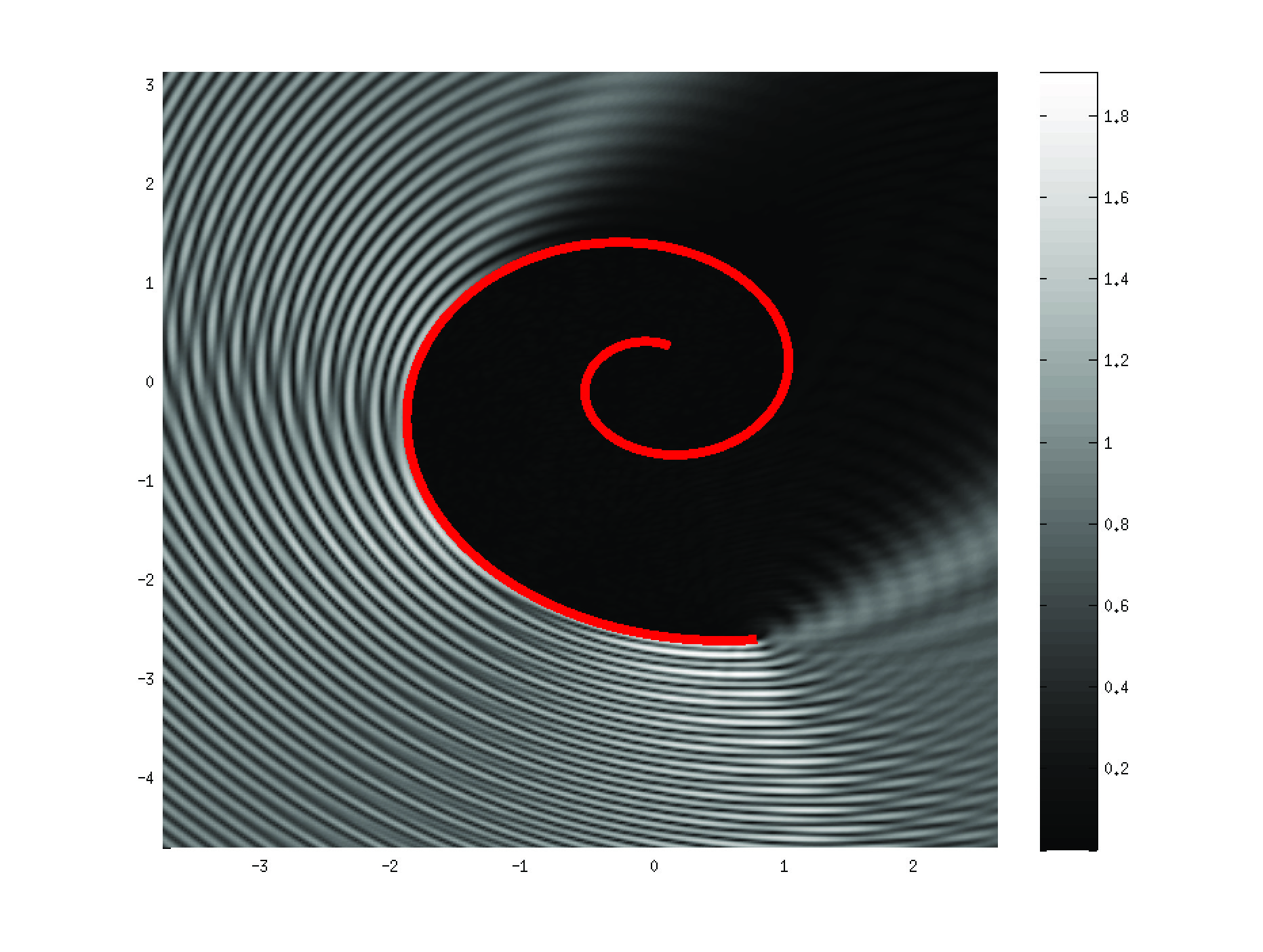} &
\includegraphics[scale=0.07]{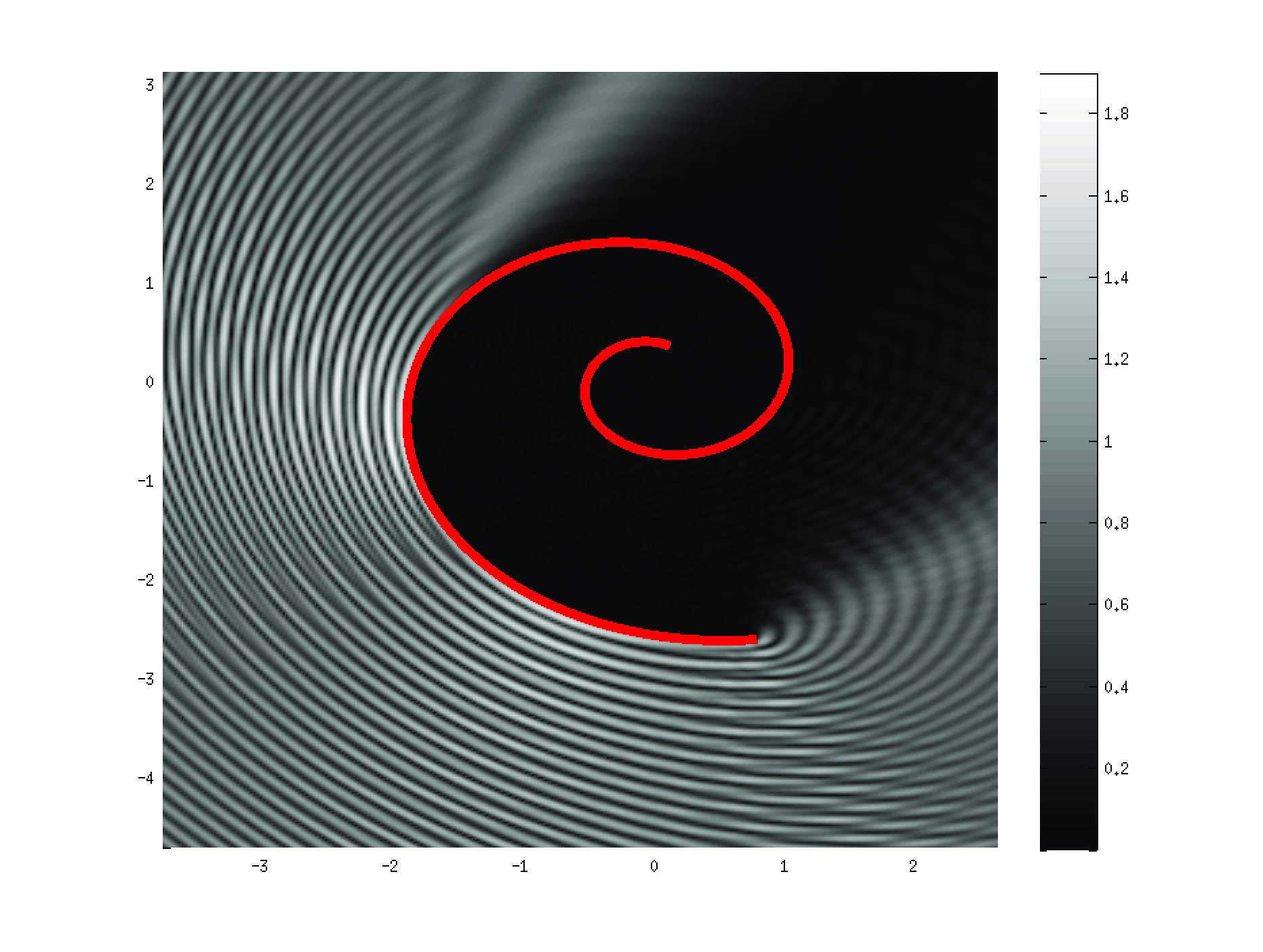} \\
(a) $U_1$ for Dirichlet Case & (b) $U_2$ for Dirichlet Case \\
\includegraphics[scale=0.07]{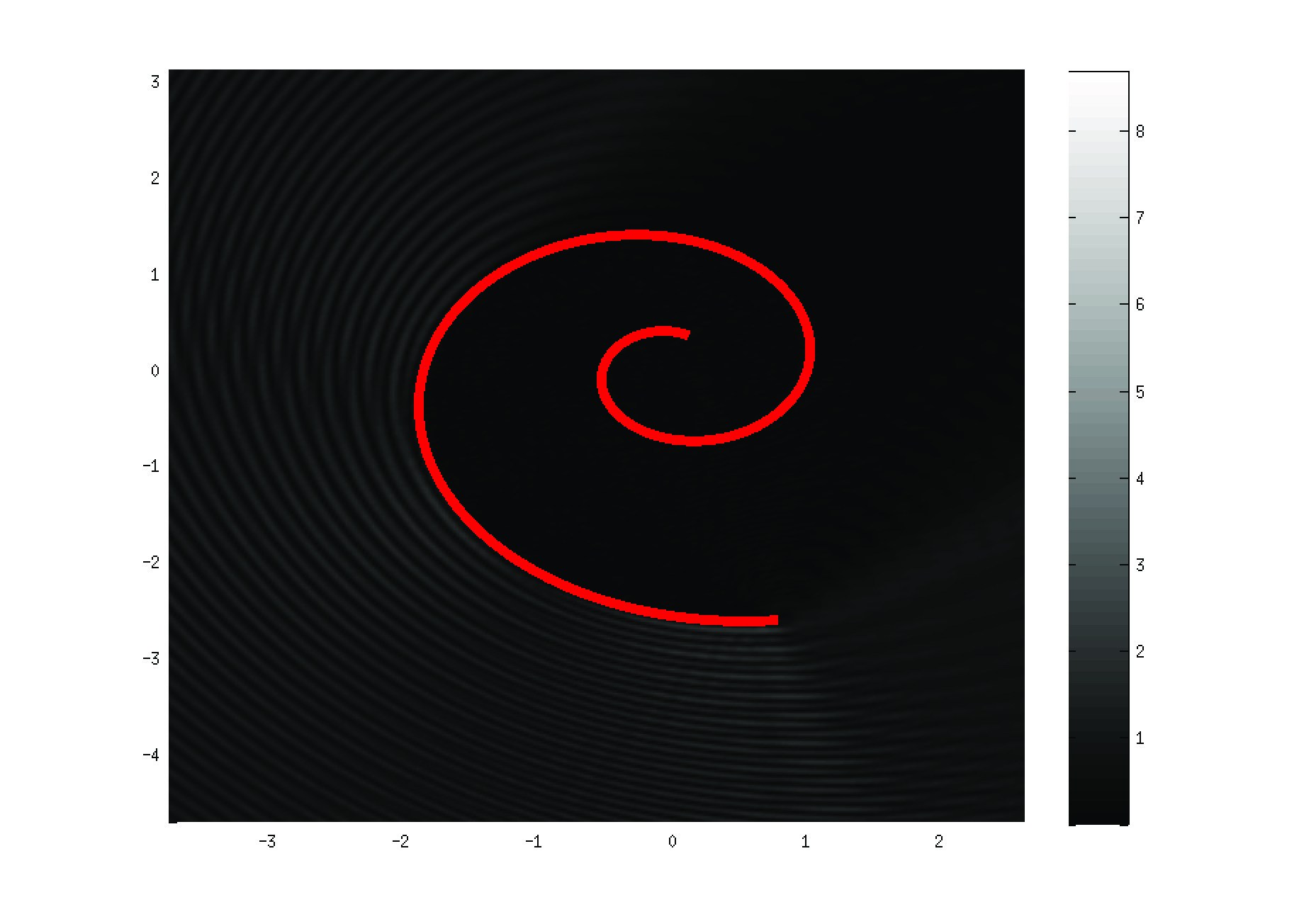} &
\includegraphics[scale=0.07]{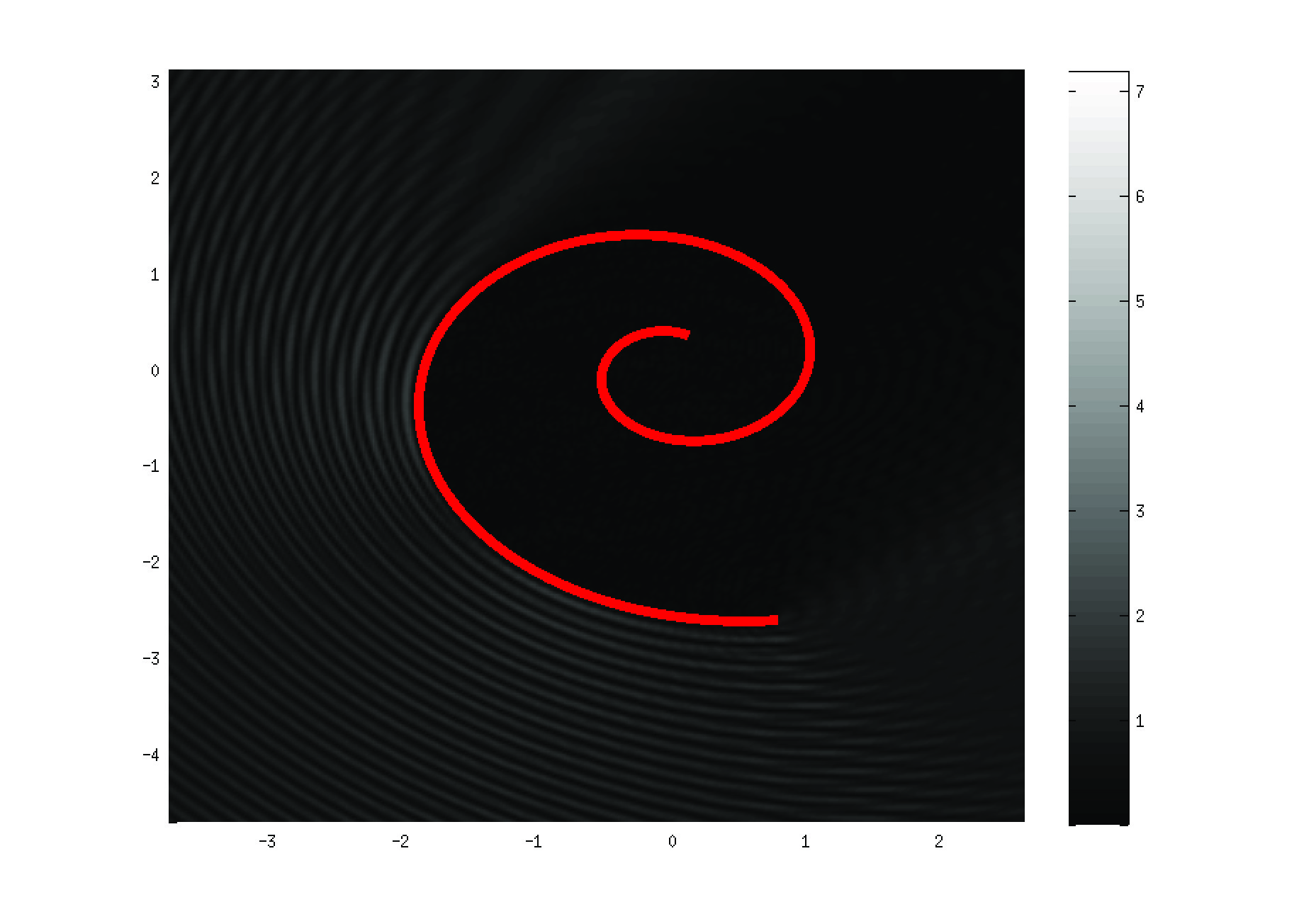} \\
(c) $U_1$ for Neumann Case & (d) $U_2$ for Neumann Case
\end{tabular}
\caption{Total field absolute values for the scattering problems by a spiral-shaped arc (Example 2).}
\label{fig:espiralfields}
\end{figure}

{\bf Example 2.} Next, we consider two single arcs given by more challenging parametrizations: a semicircle $\Gamma=\{\bs x=(\cos \frac{\pi}{2}(t+1) , \sin \frac{\pi}{2}(t+1))\in\R^2: t \in[-1,1]\}$; and, a spiral $\Gamma=\{\bs x=e^{t}(\cos 5t , \sin 5t)\in\R^2: t \in[-1,1]\}$ illuminated by the $p-$plane incident wave with incidence angle $\alpha = \frac{\pi}{2}$ and $\alpha = \frac{\pi}{4}$, respectively. The numerical error convergence for an increasing number of polynomials basis for both the Dirichlet and Neumann problems is presented in Figure~\ref{fig:errorscircle} while the corresponding total fields are plotted in Figures~\ref{fig:scirclefields} and \ref{fig:espiralfields}.

We can see from the convergence results shown in Figures~\ref{fig:errorline} and \ref{fig:errorscircle} that the proposed method can achieve more than 10 digits of accuracy. Moreover, we also infer that, after a pre-asymptotic part---depending on $N_0$ from Corollary \ref{cor:convergencerate}, and also the oscillatory behavior of the solution---, the logarithm of the error decays at a constant rate with respect to the number of polynomials. Thus, the convergence is exponential as it was stated in Corollary \ref{cor:convergencerate}.  We can also compare the results with the ones presented in \cite{BXY21}, we see that the convergence rate seems similar but we are able to achieve smaller errors with less degrees of freedom. This is more notorious for the Neumann problem where the energy norm is stronger than the uniform norm used in \cite{BXY21}. We remark however that the Nystr\"om discretization used in \cite{BXY21} should in practice be less computationally expensive than our spectral method.

\begin{figure}[htbp]
\centering
\begin{tabular}{cc}
\includegraphics[scale=0.07]{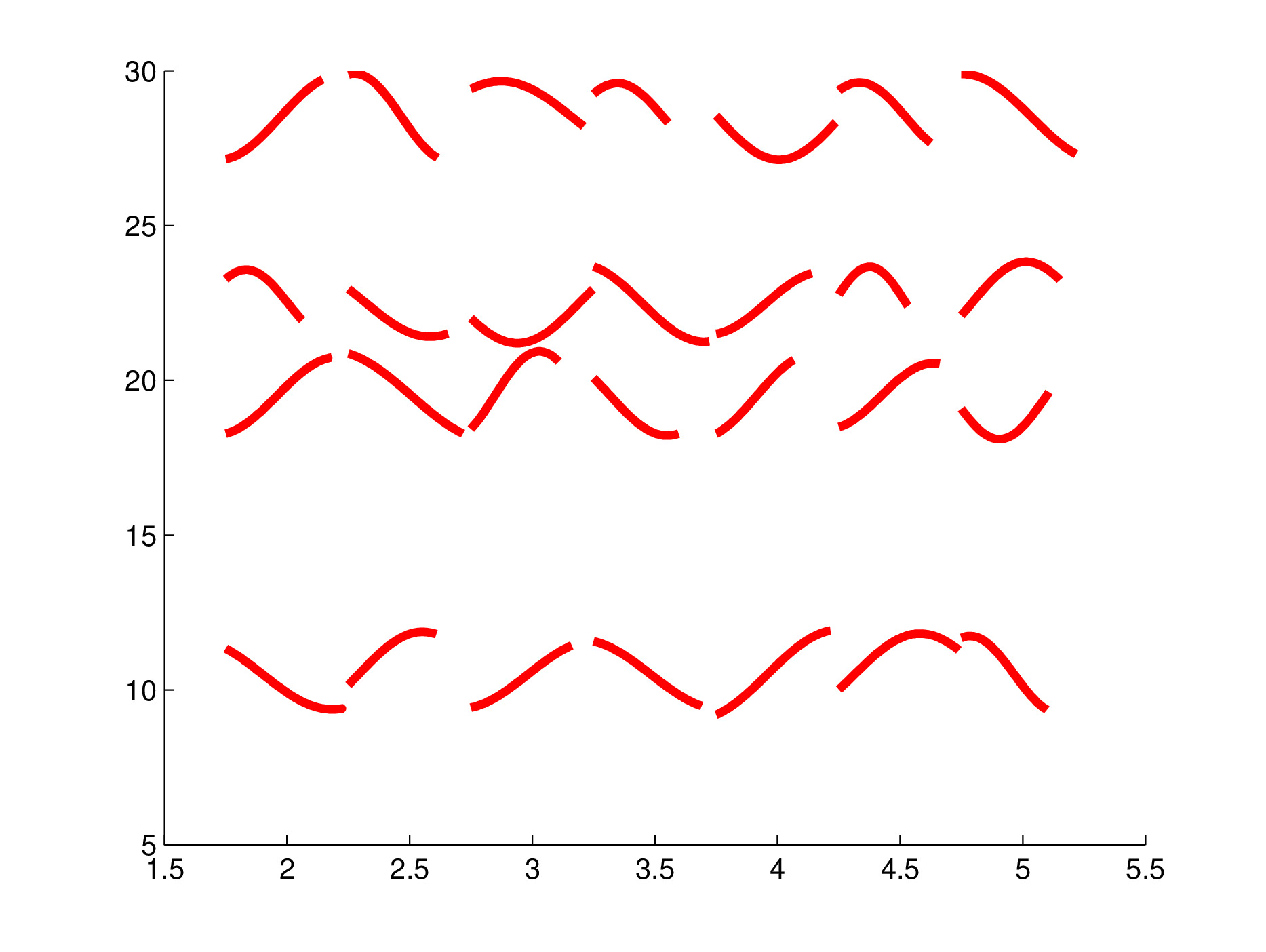} &
\includegraphics[scale=0.07]{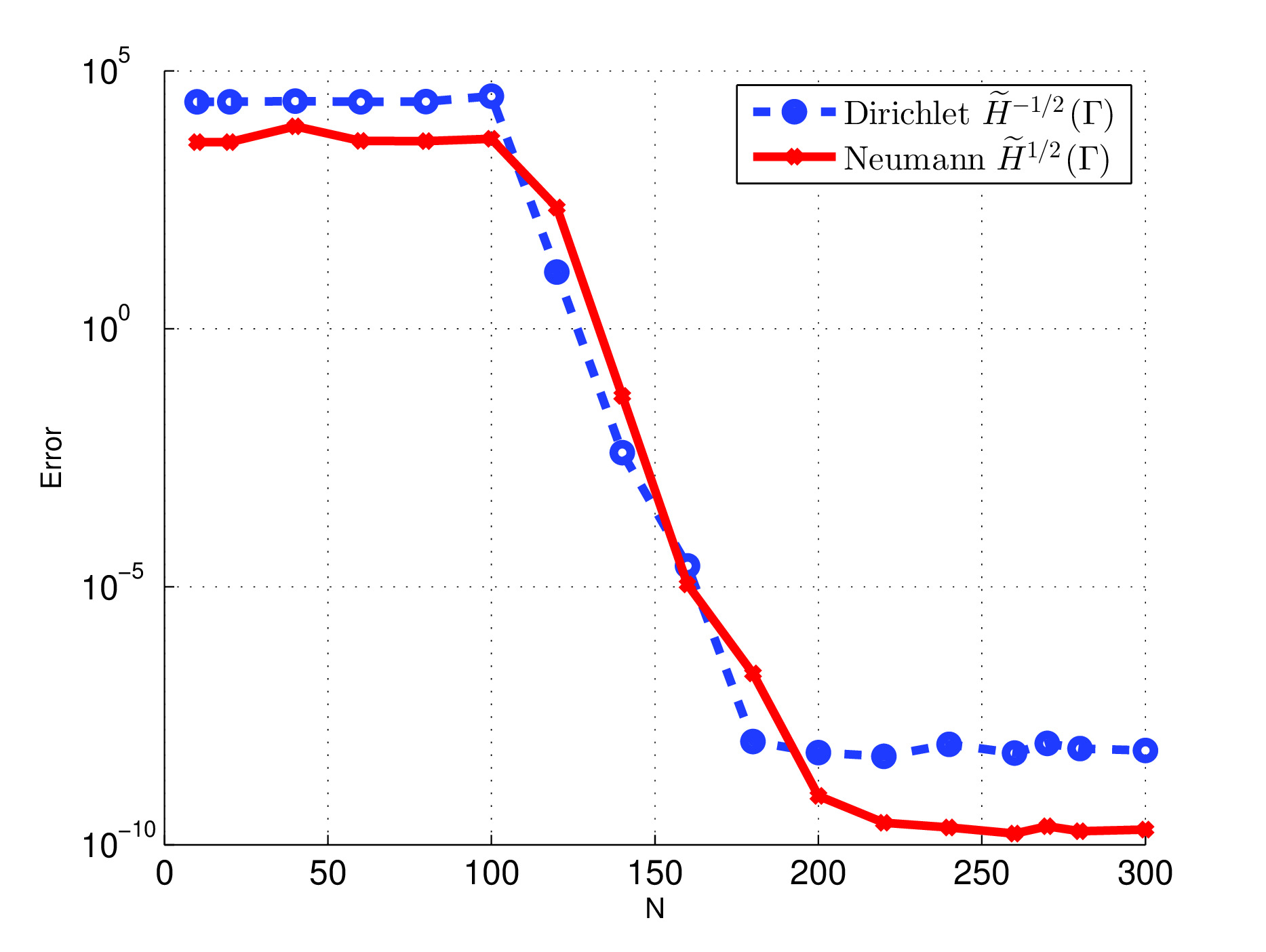} \\
(a)  &  (b)
\end{tabular}
\caption{Numerical errors (b) for the problems of scattering by 28 open arcs (a) as described in Example 3.}
\label{fig:28arcs}
\end{figure}

{\bf Example 3.} Let us now consider a more complex geometry consisting of 28 open arcs given by the general formula:
$$x(t) = at+b, \quad y(t) = c\sin(\alpha t+\gamma)+d,$$ where the real constants $a,b,c,d,\beta,\gamma$ are different for each arc, and where selected randomly on adequate ranges (see Figure~\ref{fig:28arcs}(a)). We fix an incidence angle of $\alpha =0$. In this case, we use the compression of the cross interaction matrices with {\texttt{tol}}=$10^{-10}$. We present the convergence of the numerical errors for an increasing number of polynomials basis for both the Dirichlet and Neumann problems in Figure \ref{fig:28arcs}(b). In particular, the running times are $4$ min for the Dirichlet case and $8$ min for the Neumann one, both for $N= 320$, for which the number of degrees of freedom equals to 17,920.

{\bf Example 4.}
As a fourth example, we consider the first 10 arcs of the geometry of Example 3 sorted from bottom to top and left to right in Figure \ref{fig:28arcs}, again with $\alpha  =0$, and solve the Dirichlet Problem for various values of $\omega$. As in the previous case, we use the compression algorithm with \texttt{tol}=$10^{-10}$. The results are reported in Table \ref{tab:ex4}, where $N$ denotes the polynomial degree used per arc, the error is computed in the energy norm, and NNZ\%, denotes the percentage of the matrix with non-zero entries. As expected, we observe an increase in degrees of freedom and computation times as the frequency increases. Similarly for the compression algorithm as the resolving number of $N_0$ also increases.

\begin{table}[]
\centering
\begin{tabular}{l|llll}
$\omega$ & N   & Error & NNz\% & Time(s) \\ \hline
10       & 170 & 1e-11 & 11    & 6.5     \\ \hline
50       & 240 & 1e-10 & 22    & 17      \\ \hline
100      & 310 & 1e-10 & 36    & 37      \\ \hline
150      & 400 & 1e-10 & 41    & 71      \\ \hline
208      & 520 & 1e-10 & 42    & 154     \\ \hline
250      & 610 & 1e-9  & 42    & 242
\end{tabular}
\caption{Results for a range of frequencies, errors computed against an overkill solution with a polynomial degree equal to $N+60$.}
\label{tab:ex4}
\end{table}

{\bf Example 5.}
For the last example we consider again a geometry of open arcs with the general formula used in the Example 3, $\alpha  = 0$, and \texttt{tol}=$10^{-10}$. We consider an increasing number of arcs and fix the polynomial degree per arc as $N = 200$. Results for the Dirichlet problem are reported in Table \ref{tab:ex5}, for which we observe an increase in computation times, as it should be expected.

\begin{table}[]
\centering
\begin{tabular}{l|lll}
\# Arcs & Error & NNz\% & Time(s) \\ \hline
5       & 1e-10 & 30    & 4       \\ \hline
10      & 1e-10 & 31    & 11      \\ \hline
15      & 1e-10 & 32    & 24      \\ \hline
20      & 1e-9  & 33    & 38      \\ \hline
30      & 1e-10 & 33    & 86      \\ \hline
40      & 1e-9  & 33    & 155
\end{tabular}
\caption{Results for a increasing number of open arcs, errors computed against an overkill solution with a polynomial degree equal to $260$.}
\label{tab:ex5}
\end{table}

\section{Conclusions and Future Work}
\label{sec:5}

We have presented a fast spectral Galerkin method for solving the weakly- and hyper-singular BIEs that reformulate the two-dimensional Dirichlet and Neumann problems of elastic time-harmonic scattering by multiple disjoint cracks, respectively. The numerical discretization utilizes weighted Chebyshev polynomials to treat the singular behavior of the solutions' edge singularities explicitly and, by assuming analyticity of sources and arcs geometries, exponential convergence of the numerical scheme is shown. Several numerical examples are presented to verify our theoretical results and show the accuracy and efficiency of the proposed method. Although the Nystr\"om discretization used in \cite{BXY21} displays similar numerical convergence rates, it  would be quite interesting to prove the convergence rate of the Nystr\"om method while,  considering the singular behavior of the solutions' edge singularities explicitly. Additionally, the study of appropriate spectral Galerkin method for the three-dimensional elastic problems of cracks and the application of the numerical methods for inverse problems and uncertain quantification problems will be left for future works.

\section*{Acknowledgements}
TY gratefully acknowleges support from NSFC through Grant No. 12171465.

\bibliographystyle{siam}
\bibliography{references}

\end{document}